\documentclass[twocolumn, 10pt, journal,letterpaper]{IEEEtran}
\usepackage{url}
\usepackage{graphicx}
\usepackage{geometry}
\usepackage{nomencl}
\usepackage{algorithm}
\usepackage[noend]{algpseudocode}
\usepackage{etex}
\usepackage{tikz}
\usetikzlibrary{shapes,arrows}
\usepackage{pstricks}
\usepackage{pst-node,pst-blur}
\usetikzlibrary{arrows}
\usepackage{amsfonts}
\usepackage{tabularx}
\usepackage{subfigure}
\usepackage{amssymb}
\usepackage{amsmath}
\usepackage{booktabs}
\usepackage[font=footnotesize,skip=1pt]{caption}
\usepackage{multirow}
\usepackage{multicol}
\usepackage{steinmetz}
\usepackage{epsfig}
\usepackage{epstopdf}
\usepackage{ntheorem}
\usepackage{xspace}
\usepackage{setspace}
\usepackage{breqn}

\usepackage{cite}
\newtheorem{theorem}{Theorem}
\newtheorem{remark}{Remark}

\newtheorem{lemma}{Lemma}

\newtheorem{corollary}{Corollary}
\newtheorem*{proof}{Proof}

\usepackage{stfloats} 
\usepackage[cmex10]{mathtools}
\usepackage{physics}

\allowdisplaybreaks
\makeatletter
\def\@eqnnum{{\normalsize \normalcolor (\theequation)}}
\makeatother
\begin{document}
\title{Optimal Reflection Coefficients for ASK Modulated Backscattering from Passive Tags}
\author{Amus Chee Yuen Goay,
\IEEEmembership{Student Member, IEEE}, Deepak Mishra, \IEEEmembership{Senior Member, IEEE}, Aruna Seneviratne, \IEEEmembership{Senior Member, IEEE}
\thanks{A. C. Y. Goay, D. Mishra, and A. Seneviratne are with the School of Electrical Engineering and Telecommunications, University of New South Wales, Sydney, NSW 2052, Australia (email: a.goay@student.unsw.edu.au, d.mishra@unsw.edu.au, a.seneviratne@unsw.edu.au). }
\thanks{A preliminary conference version of this work was presented at IEEE MASCOTS Workshop in Nice, France, in October 2022~\cite{goay2022ber}.}\vspace{-0.5mm}}

\maketitle

\begin{abstract}
This paper studies backscatter communication (BackCom) systems with a passive backscatter tag. The effectiveness of these tags is limited by the amount of energy they can harness from incident radio signals, which are used to backscatter information through the modulation of reflections. To address this limitation, we adopt a practical Constant-Linear-Constant (CLC) energy harvesting model that accounts for the harvester's sensitivity and saturation threshold, both of which depend on the input power. This paper aims to maximize this harvested power at a passive tag by optimally designing the underlying M-ary amplitude-shift keying (ASK) modulator in a monostatic BackCom system. Specifically, we derive the closed-form expression for the global optimal reflection coefficients that maximize the tag’s harvested power while satisfying the minimum symbol error rate (SER) requirement, tag sensitivity, and reader sensitivity constraints. We also proposed optimal binary-ASK modulation design to gain novel design insights on practical BackCom systems with readers having superior sensitivity. We have validated these nontrivial analytical claims via extensive simulations. The numerical results provide insight into the impact of the transmit symbol probability, tag sensitivity constraint, and SER on the maximum average harvested power. Remarkably, our design achieves an overall gain of around $13\%$ over the benchmark, signifying its utility in improving the efficiency of BackCom systems. Moreover, our proposed solution methodology for determining the maximum average harvested power is applicable to any type of energy harvesting model that exhibits a monotonic increasing relationship with the input power.
\end{abstract}


\vspace{-0.5mm}

\section{Introduction}
\bstctlcite{IEEEexample:BSTcontrol}

\vspace{-0.5mm}

\IEEEPARstart{B}{\lowercase{ackscatter}} communication (BackCom) is a low-power, cost-effective wireless technology that enables tags to transmit data by reflecting incident radio frequency (RF) signals to a receiver. The concept of backscattering was first applied during World War II to distinguish between allies' and enemies' aircraft. In the 1960s, it found commercial use, with Radio Frequency Identification (RFID) becoming one of its most widespread applications. More recently, BackCom has gained significant attention as a green solution for low-power communication systems~\cite{niu2019overview}.\footnote{Generally, the BackCom systems can be classified based on their read range and system configuration. These categories include near-field and far-field systems, which differ in both operating frequency and hardware design. Near-field BackCom, typically operating at $125$ kHz or $13.56$ MHz, relies on magnetic coupling for information transfer. In contrast, far-field BackCom operates in the $860–960$ MHz frequency range and uses electromagnetic radiation to transmit data over longer distances.}

BackCom is playing a key role in the growth of the Internet of Things (IoT) due to its ability to drastically reduce both manufacturing and operational costs for IoT devices~\cite{landaluce2020review}, while also extending their operational lifespan. A recent survey on BackCom-based green IoT systems highlights its potential for joint sensing and wireless communication, providing insights into the technology's working principles, applications, and associated challenges~\cite{toro2021backscatter}. Despite its promising applications, the low efficiency of BackCom in far-field scenarios remains a major bottleneck~\cite{kaur2011rfid,want2006introduction}. This bottleneck restricts the performance of BackCom systems, as the energy available to power the passive tags is often insufficient. To overcome this challenge, careful selection of design parameters is required to enhance energy harvesting by the tags and meet the communication requirements of BackCom systems.

\subsection{State-of-the-Art}\label{sec:RW}

The performance of passive BackCom systems is typically evaluated based on key metrics such as data rate, tag-to-reader transmission range, power transfer, and bit error rate (BER)~\cite{nikitin2006performance,muralter2019selecting}. The data rate is influenced by the symbol rate and the number of bits per symbol, while the transmission range is determined by factors like tag sensitivity, reader sensitivity, and BER. Optimizing these metrics is crucial for improving the overall efficiency and reliability of BackCom systems. Unlike conventional wireless communication devices, passive backscatter tags lack active RF components and instead utilize load modulation to transmit information.\footnote{Backscatter tags can be classified into three main categories: passive, semi-passive, and active\cite{lu2018ambient}. Passive tags are batteryless and rely on external energy sources to power their circuits. On the other hand, semi-passive and active tags are equipped with internal batteries. While semi-passive tags still use the incident RF signal for backscattering, active tags can generate their own signals independently.} The tag modulates the backscattered signal by switching the load impedance at the antenna terminal~\cite{goay2022ask}. The choice of load impedance is directly related to the modulation scheme employed. The three most common modulation techniques used in BackCom systems are amplitude-shift keying (ASK), phase-shift keying (PSK), and frequency-shift keying (FSK)~\cite{de2005design,fasarakis2015coherent}.

One of the primary advantages of load modulation is its low power consumption and simple circuit design. However, the tag's performance is highly dependent on the selection of the load impedance. Studies have shown that this load selection significantly impacts the overall performance of BackCom systems~\cite{rao2005impedance,nikitin2005power,nikitin2007differential}. Bletsas et al.~\cite{bletsas2010improving} proposed a load selection strategy to minimize BER for ASK and PSK modulations without considering tag power sensitivity. On the other hand, De Vita et al.~\cite{de2005design} introduced a method for selecting load impedance with equal impedance mismatch in both states. Meanwhile, Karthaus et al.~\cite{karthaus2003fully} examined load impedance selections explored in~\cite{nikitin2005power,bletsas2010improving,de2005design} and showed that the harvested power is directly influenced by modulation depth. In recent studies, reflection coefficient optimization has emerged as a key focus for improving performance across different modulation schemes. In~\cite{goay2023qos,goay2024qos}, we demonstrated optimal reflection coefficient selection that accounts for BackCom system operational requirements, enhancing performance for both PSK and quadrature amplitude modulation (QAM) in various applications. Additionally, a novel selective FSK modulation is proposed to maximize transmission range by determining the optimal reflection coefficient~\cite{goay2024backscatter}.

Moreover, there is increasing interest in enhancing the data rates of backscatter devices~\cite{khadka2022index}. Specifically, achieving higher data rates through increasing the number of bits per symbol is often preferred over increasing the baseband frequency, as the latter leads to higher power consumption. Recent research has explored the use of high-order modulation techniques to increase data rates. For instance, ambient BackCom systems employing M-ary phase-shift keying (M-PSK) and M-ary frequency-shift keying (M-FSK) modulation schemes have been studied to enhance data transmission rates~\cite{qian2018iot,tao2019ambient}. Additionally, the selection of reflection coefficients for high-order QAM modulation at the tag has been proposed to ensure equal Euclidean distances in the constellation diagram~\cite{thomas2010qam}. Another example is the 4-QAM modulator for RFID tags designed to minimize power loss while maintaining low BER~\cite{boyer2012coded}. In addition to these modulation techniques, spatial modulation has been explored in ambient BackCom systems to achieve both high data rates and improved spectrum efficiency~\cite{niu2019spatial,raghavendra2022generalized}.

Furthermore, various methods have been explored to enhance BackCom system performance. For instance, the work in~\cite{kimionis2014increased} investigated a bistatic architecture to maximize the read range. A novel time-division multiple access protocol for a tag-to-tag cooperative scheme was proposed in~\cite{10118930} to improve throughput in a 3-tag BackCom system. Kimionis et al.~\cite{kimionis2014enhancement} suggested using reflection amplifiers to improve the signal-to-noise ratio (SNR). In~\cite{6685977}, unequal forward error correction coding was explored to balance the trade-off between harvested power and BER. A joint optimization of the reflection coefficient for maximizing system throughput was studied in~\cite{Mishra_TCOM_Multitag}.

\subsection{Motivation and Contributions}\label{sec:motiv}

The BackCom system performs poorly in far-field applications because the harvested energy at the tag decreases dramatically over longer distances~\cite{finkenzeller2010rfid}. Thus, the passive tag's transmission range is generally short and cannot execute advanced tasks. Therefore, the utility of the tag can be significantly improved by maximizing the harvested power at the tag. This will also enable the tag to perform more onboard operations and support more applications. Our objective is to develop a high-data-rate BackCom system using the high-order M-ary ASK ($M$-ASK) modulation scheme. Unlike existing works that consider equal probability for the transmit symbols during backscattering, our investigation focuses on maximizing the average harvested power with unequal symbol probabilities.

The studies in~\cite{nikitin2005power,nikitin2007differential,bletsas2010improving,de2005design,karthaus2003fully} have explored various load selections without finding the optimal value for enhancing the tag performance. \textit{To the best of our knowledge, maximizing the tag's average harvested power through optimal reflection coefficient selection for $M$-ASK modulation, under symbol error rate (SER), reader sensitivity, and energy constraints, has not yet been investigated.} The main contributions of this article are as follows:
\begin{enumerate}
    \item We detail the load-dependent performance metrics, including harvested power, backscattered power, and SER. We then formulate the problem of maximizing the backscatter tag's harvested power by jointly optimizing the reflection coefficients for the $M$-ASK modulation scheme, while accounting for practical operational requirements.
    \item To tackle the non-convex nature of the problem, we decouple it into two parts by converting an inequality constraint into an equality constraint, thereby reformulating the problem into a convex one. We then iteratively solve for the optimal solution, determining all solution sets corresponding to different transmit symbol probability sequences.
    \item A specific case of the BackCom system utilizing binary ASK (BASK) modulation, tailored for RFID applications, is investigated. We derived a closed-form solution by applying the Karush-Kuhn-Tucker (KKT) conditions. Given its simplicity, this BASK modulation design is well-suited for low-cost sensing and low-data-rate communication systems.
    \item Simulation results are presented to quantify the maximum average harvested power for various applications, under different values of the key system parameters.Design insights regarding optimal reflection coefficients are provided, validating the effectiveness of the proposed optimal solution.
    \item We analyze the saturation nonlinear energy harvesting model and compare it with the adopted Constant-Linear-Constant (CLC) model to assess performance differences. Our results show that the CLC model serves as a practical approximation for energy harvesting systems. Moreover, our proposed algorithm can be applied to any energy harvesting model that exhibits a monotonically increasing behavior with input power.
\end{enumerate}

\subsection{Paper Organization and Notations}

This paper is organized as follows. Section~\ref{sectionSD} describes the BackCom system model and the signal transmission process of the backscatter tag, along with its modulation technique. In Section~\ref{sectionPMB}, the performance metrics of backscatter tags are discussed in more detail. Sections~\ref{sectionPD} and~\ref{Solution} outline the optimization problems and the proposed solution methodology, respectively. Section~\ref{sectionRD} presents the numerical results, followed by the conclusions in Section~\ref{sec:conclusion}

Notations: Throughout this paper, $\Re \left( \cdot \right)$ and $\Im \left( \cdot \right)$ denote the real and imaginary parts of a complex number, respectively. Bold upper-case letter represent matrix, while lower-case letter with subscript $nm$ refer to the ($n,m$)-th entry of a matrix. The operator $\mathbb{E}\left[\cdot\right]$ denotes statistical expectation. Additionally, $\abs{\cdot}$ returns the absolute value, and $j = \sqrt{-1}$ represents the imaginary unit. 

\vspace{-1mm}

\section{System Description}\label{sectionSD}

\subsection{System Model and Transmission Protocol}

\begin{figure}[t!]
	\centering
\includegraphics[width=3.3in]{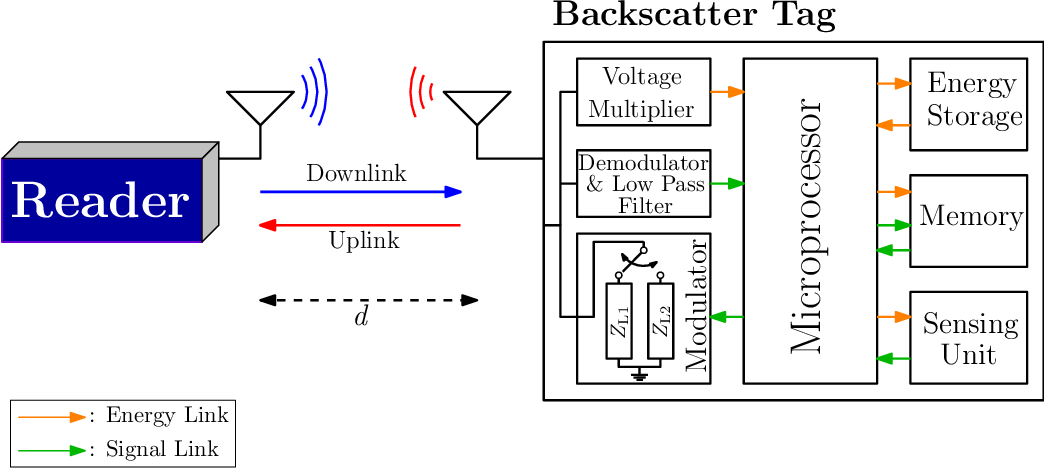}
	\caption{Monostatic BackCom system model.}
	\label{fig1}
\vspace{-2mm}
\end{figure}

This paper studies a monostatic BackCom system with one reader and one passive tag separated by distance $d$, as shown in Fig.~\ref{fig1}. As a dedicated power source, the reader stably broadcasts an unmodulated RF carrier with constant power to the passive tag in the downlink. Then, the passive backscatter tag\footnote{There are various passive backscatter tags, of which the RFID tag is the most well-known one. In the latest research, the traditional RFID tag integrated with sensing electronics, transforming it into a sensing and computational platform, has been studied for IoT applications. The tag with sensing capability is called computational RFID (CRFID), which has higher power consumption during operation~\cite{landaluce2020review}.}, functioning as an IoT device with sensing capability, modulates the backscattered signal with its unique identifier and the sensing data. This modulated signal is then transmitted to the reader in the uplink. As depicted in Fig.\ref{fig1}, the passive tag comprises an antenna, voltage multiplier (VM), demodulator, low-pass filter (LPF), modulator, microprocessor, energy storage, memory, and a sensing unit. 

When a sinusoidal electromagnetic (EM) wave is impinging the tag antenna, the VM rectifies the induced voltage into DC power and delivers it to the microprocessor. The induced voltage is used to generate the low-frequency subcarrier, which is then used to produce the modulated subcarrier with baseband signal~\cite{finkenzeller2010rfid}. Once the microprocessor is powered, the modulator encodes the backscattered signal by switching between different load impedances based on system requirements and the modulation scheme. This paper focuses on a BackCom system utilizing ASK modulation. The system model can be extended to support multiple tags with minimal modifications by incorporating a medium access control $\left(\text{MAC}\right)$ protocol. In most BackCom systems, the ALOHA protocol is generally used as the MAC for anti-collision~\cite{wu2013binary,tegos2020slotted}.

\vspace{-1mm}

\subsection{M-ary ASK Modulation}
\begin{figure}[t!]
    \centering
    \includegraphics[width=2.4in]{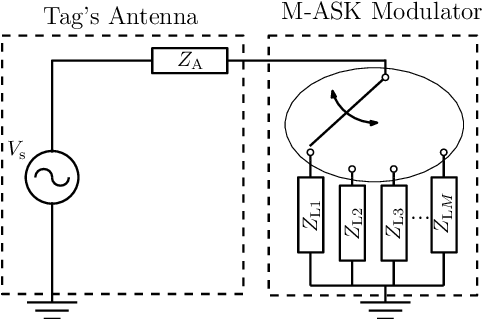}
    \caption{Equivalent Thevenin Circuit of the Backscatter Tag with Minimum Scattering Antenna for $M$-ASK Modulation.}
    \label{fig12abc}
    \vspace{-2mm}
\end{figure}

This section introduces the $M$-ASK modulation scheme for the backscatter tag, where $M = 2^n$ and $n \in \{1, 2, 3, 4, \dots \}$ represents the number of bits per transmit symbol. For example, when $n=1$, the transmit symbols are \{‘1’, ‘0’\}, and when $n=2$, the symbols are \{‘11’, ‘10’, ‘01’, ‘00’\}. Since the tag utilizes load modulation, it switches between different load impedances to transmit each symbol. Thus, there are $M$ distinct load impedances $\left(Z_{\rm L1}, Z_{\rm L2}, Z_{\rm L3}, \dots, Z_{{\rm L}M}\right)$, each corresponding to a specific symbol $\left(S_1, S_2, S_3, \dots, S_M\right)$ in the $M$-ASK modulation scheme. An equivalent Thevenin circuit of the backscatter tag using $M$-ASK modulation is shown in Fig.~\ref{fig12abc}. Without loss of generality, we assign the load impedances $Z_{{\rm L}i}$ to represent symbols in descending order of their decimal values. For instance, in 4-ASK, the symbols \{‘11’, ‘10’, ‘01’, ‘00’\} correspond to the load impedances $Z_{\rm L1}, Z_{\rm L2}, Z_{\rm L3}$, and $Z_{\rm L4}$, respectively.

\vspace{-1mm}
 
\section{Performance Metrics for BackCom}\label{sectionPMB}

The tag sensitivity, reader sensitivity, and SER are identified as key factors determining the maximum transmission range of BackCom systems~\cite{nikitin2006performance, muralter2019selecting}. In this section, we consider these factors and incorporate them into our analysis of reflection coefficient selection.

\vspace{-1mm}

\subsection{Tag Power Sensitivity}\label{sectionMLPTATIC}

The backscatter tag modulates its backscattered signal by altering its circuit impedance, and the amount of energy harvested by the tag depends on the load impedance at each state. The reflection coefficient is commonly used to characterize the tag's circuit impedance during backscattering. According to Kurokawa~\cite{kurokawa1965power}, the reflection coefficient is the ratio of the reflected power wave to the total incident power wave. We denote the reflection coefficient as $\Gamma_i$, where $i \in \mathbb{M} = \{1,2,3,...,M\}$, which corresponds to the tag connecting to the load impedance $Z_{{\rm L}i}$ and transmit symbol $S_i$. The reflection coefficient $\Gamma_i$ is expressed as~\cite{kurokawa1965power}:
\begin{equation}
    \Gamma_i \triangleq \frac{Z_{{\rm L}i} - \bar{Z}_{\rm A}}{Z_{{\rm L}i} + Z_{\rm A}}, \quad \forall i \in \mathbb{M},
    \label{eq1}
\end{equation}
 where $Z_{\rm A} = R_{\rm A} + j X_{\rm A}$ is the antenna impedance, $\bar{Z}_{\rm A}$ is the conjugate of $Z_{\rm A}$, and $Z_{{\rm L}i} = R_{{\rm L}i} + j X_{{\rm L}i}$. Here, $R_{{\rm L}i}$ and $R_{\rm A}$ are the load resistance and antenna resistance, respectively, while $X_{{\rm L}i}$ and $X_{\rm A}$ are the load reactance and antenna reactance. To simplify the design optimization, we substitute the expression of $Z_{\rm A}$, $\bar{Z}_{\rm A}$ and $Z_{{\rm L}i}$ into~\eqref{eq1} and rewrite it in the rectangular form as: 
 \begin{align}
\Gamma_i &= \frac{R^2_{{\rm L}i} - R^2_{\rm A} + \left(X_{{\rm L}i} + X_{\rm A}\right)^2 + j 2 R_{\rm A} \left(X_{{\rm L}i} + X_{\rm A}\right)}{\left(R_{{\rm L}i} + R_{\rm A}\right)^2 + \left(X_{{\rm L}i} + X_{\rm A}\right)^2} \nonumber \\
        &= \Gamma_{{\rm a}i} + j \Gamma_{{\rm b}i}, \quad \forall i \in \mathbb{M},
    \label{eq3}
\end{align}
where $\Gamma_{{\rm a}i} = \frac{R^2_{{\rm L}i} - R^2_{\rm A} + \left(X_{{\rm L}i} + X_{\rm A}\right)^2}{\left(R_{{\rm L}i} + R_{\rm A}\right)^2 + \left(X_{{\rm L}i} + X_{\rm A}\right)^2}$ and $\Gamma_{{\rm b}i} = \frac{ 2 R_{\rm A} \left(X_{{\rm L}i} + X_{\rm A}\right)}{\left(R_{{\rm L}i} + R_{\rm A}\right)^2 + \left(X_{{\rm L}i} + X_{\rm A}\right)^2}$. The reflection coefficient plays a crucial role in determining the tag's performance in a BackCom system. Knowing that the tag can harvest energy and backscatter signal simultaneously, the power allocation during backscattering depends on $\Gamma_i$. We denote $P_{{\rm T}i}$ as the input power at the tag, which is given by~\cite{goay2022ask}:
\begin{align}
    P_{{\rm T}i} &\triangleq P_{\rm a} \left(1 - \abs{\Gamma_i}^2\right) \nonumber \\
            & \hspace{-1mm} \stackrel{(s_1)}{=} P_{\rm a} \left(1-\Gamma^2_{{\rm a}i}-\Gamma^2_{{\rm b}i}\right), \quad \forall i \in \mathbb{M},
    \label{eq4}
\end{align}
where $(s_1)$ is obtained using \eqref{eq3} and $P_{\rm a}$ is the maximum available power. We use the close-in free space reference distance model\footnote{This model can be applied to the 0.5-100 GHz frequency band~\cite{haneda20165g} and can be reduced to Friis transmission model (free space path loss) when we set $n=2$.}~\cite{andersen1995propagation,feuerstein1994path}, giving $P_{\rm a}$ as:
\begin{equation}
    P_{\rm a} \triangleq P_{\rm t} G_{\rm t} G_{\rm r} \left( \frac{\lambda}{4 \pi d_{\rm o}} \right)^2 \left( \frac{d_{\rm o}}{d} \right)^n
\end{equation}
where $d_{\rm o} = 1$m is the reference distance~\cite{maccartney2017rural}, $P_{\rm t}$ is the transmission power of the reader's antenna, $G_{\rm t}$ and $G_{\rm r}$ are the antenna gains of the tag and reader, $n$ is the path loss exponent, and $\lambda = \frac{c}{f}$ is the wavelength of the RF carrier with $c$ as the speed of light and $f$ is the carrier frequency.

In practical scenarios, the tag will not harvest any power if the input power is below a certain threshold, called the tag's power harvesting sensitivity $P_{\rm sen}$. This sensitivity arises from the minimum voltage required by the diodes in the rectifier to operate~\cite{valenta2014harvesting}. When the input power exceeds a certain level, the harvested power becomes constant, reaching the tag’s saturation point. This upper threshold, referred to as the saturated harvested power $P_{\rm sat}$, beyond which no further power can be harvested.

To ensure reliable performance and efficient power usage in BackCom systems, it is essential to account for both $P_{\rm sen}$ and $P_{\rm sat}$. We adopt a Constant-Linear-Constant (CLC) energy harvesting model, which defines the harvested power $P_{{\rm L}i}$ as a function of the input power, characterizing the system's behavior in different power regions~\cite{alevizos2018sensitive}:
\begin{equation}
    P_{{\rm L}i} = 
    \begin{cases}
        0, & P_{{\rm T}i} \in \left[ 0, P_{\rm sen} \right], \\
        \eta_{\rm h} \left( P_{{\rm T}i} - P_{\rm sen} \right), & P_{{\rm T}i} \in \left[ P_{\rm sen}, P_{\rm sat} \right], \\
        \eta_{\rm h} \left( P_{\rm sat} - P_{\rm sen} \right), & P_{{\rm T}i} \in \left[\right. \hspace{-1mm} P_{\rm sat}, \infty \hspace{-1mm} \left.\right), 
    \end{cases}
\end{equation}
for $i \in \mathbb{M}$ and $\eta_{\rm h} \in \left[0,1\right]$ represents the tag's energy harvesting efficiency. This work focuses on improving passive BackCom system performance, addressing the major limitation of short communication range due to passive tags relying on harvested energy from the reader to power their circuitry. Therefore, the selection of the reflection coefficient is critical to ensure the tag operates effectively at the maximum specified range. Our focus is on maximizing the average harvested power when the tag is at its maximum transmission distance. At this distance, the input power received by the tag is expected to be below $P_{\rm sat}$ due to significant power attenuation with increasing distance. Therefore, in this work, we focus exclusively on scenarios where $P_{{\rm T}i} \in \left[ P_{\rm sen}, P_{\rm sat} \right]$.

In general, the operation of a passive backscatter tag's circuitry can be categorized into two models: continuous operation (non-duty-cycled operation) and the harvest-then-transmit (HTT) protocol\footnote{The WISP tag is a prime example of a passive backscatter tag operating on a duty cycle. With additional sensing, computation, and data storage capabilities, the WISP tag consumes more power than conventional passive RFID tags. Thus, it uses a duty cycle approach that alternates between energy harvesting and backscattering. The operational power cycle is shown in~\cite{sample2008design}.} (duty-cycled operation)~\cite{shi2020energy,lu2016wireless}. The first model allows the tag to simultaneously harvest energy and power its circuitry during backscattering. In contrast, the HTT model requires a separate energy harvesting period before operation. Specifically, the HTT model utilizes power from both the incident EM wave and onboard energy storage, which is charged during the energy harvesting phase. The choice between these designs depends on various factors, such as the operating environment, voltage, and power requirements of the tag. 

The tag remains off and only activates when sufficient power and a minimum threshold voltage are available. During operation, the tag needs to meet a minimum harvested power threshold $P_{\rm L,min}$ to sustain signal backscattering. In the case of continuous operation, the tag relies on $P_{\rm L,min}$ to maintain its circuitry operation. On the other hand, in the HTT model, the tag uses the incident RF carrier to generate the clock signal~\cite{finkenzeller2010rfid}, requiring $P_{\rm L,min}$ to create the baseband signal. In both models, the tag remains inactive and cannot generate any information when $P_{{\rm L}i} < P_{\rm L,min}$. Thus, ensuring that $P_{{\rm L}i} \geq P_{\rm L,min}$ is crucial for the sustainability of the BackCom system.

\subsection{Reader Sensitivity}

The second key factor in tag design is reader sensitivity, which refers to the minimum power $P_{\rm r,min}$ required for the reader to successfully decode and retrieve data from the received signal in a given noise environment. In a BackCom system, the received signal at the reader is the backscattered signal from the tag, which has propagated through the transmission channel. Thus, the power of the received signal $P_{{\rm r}i}$ at the reader is directly proportional to the power of the backscattered signal $P_{{\rm b}i}$ from the tag. The received power at the reader is given by:
\begin{equation}
    P_{{\rm r}i} \triangleq P_{{\rm b}i} G_{\rm r} \left( \frac{\lambda}{4 \pi d_o} \right)^2 \left(\frac{d_o}{d}\right)^n, \quad \forall i \in \mathbb{M},
\end{equation}

While many previous studies simplify their analyses by not considering structural mode scattering in the backscattered signal, we address this limitation in our work. Specifically, we take into account the impact of structural mode scattering and assume the tag is equipped with a minimum scattering antenna~\cite{fuschini2008analytical}.\footnote{We have introduced a signal subtraction technique that makes the minimum scattering antenna assumption practical, ensuring that the optimal reflection coefficients determined under this assumption can be applied to all types of antennas without performance degradation~\cite{goay2024tag}.} Accordingly, the power of the backscattered signal $P_{{\rm b}i}$ is given by~\cite{nikitin2006theory}:
\begin{align}
    P_{{\rm b}i} &\triangleq & \hspace{-3mm} &\eta_{\rm b} P_{\rm a} G_{\rm t} \abs{1- \Gamma_i}^2 \nonumber \\
    &\stackrel{(s_2)}{=} & \hspace{-3mm} &\eta_{\rm b} P_{\rm a} G_{\rm t} \Bigl[ \left( 1 - \Gamma_{{\rm a}i} \right)^2 + \Gamma^2_{{\rm b}i} \Bigr], \quad \forall i \in \mathbb{M},
\end{align}
where $\left(s_2\right)$ is obtained using \eqref{eq3} and $\eta_{\rm b} \in \left[ 0,1 \right]$ is the tag's backscattering efficiency. The $P_{{\rm b}i}$ is load-dependent, varying significantly with $\Gamma_i$. Therefore, reader sensitivity sets a constraint on the selection of the reflection coefficient at the tag. In this work, we focus on optimizing tag design by ensuring that $P_{{\rm b}i}$ meets the minimum required power $P_{\rm b,min}$ to satisfy reader sensitivity. Thus, $P_{{\rm b}i} \geq P_{\rm b,min}$ becomes a BackCom system operation requirement.

\subsection{Symbol Error Rate}

The third factor limiting the tag performance is the SER, defined as the number of symbols misidentified by the reader out of the total number of transmitted symbols at a given time interval~\cite{fuschini2008efficiency}. The probability $P^{\left(i,k\right)}_{\rm e}$ of symbol $S_i$ being misidentified as symbol $S_k$ is given by the following equation~\cite{fuschini2008analytical}:
\begin{equation}
    P^{\left(i,k\right)}_{\rm e} \triangleq \frac{1}{2} \text{erfc}\left( \frac{\abs{V_0} m_{i,k}}{2 \sqrt{2} \sigma} \right), \hspace{2mm} \forall i,k \in \mathbb{M}, \hspace{1mm} i \neq k,
    \label{eq5}
\end{equation}
where $V_0=\sqrt{8 R_{\rm r} P_{\rm r, m}}$ is the induced voltage at the reader's antenna when the tag is scattering EM wave at the perfect matched condition $\left(Z_{{\rm L}i} = \bar{Z}_{\rm A}\right)$~\cite{karthaus2003fully}, where $R_{\rm r}$ representing the reader's antenna resistance and $P_{\rm r, m} = \eta_{\rm b} P_{\rm a} G_{\rm t} G_{\rm r} \left(\tfrac{\lambda}{4 \pi d_o}\right)^2 \left(\tfrac{d_o}{d}\right)^n$. The inevitable additive white Gaussian noise $n_r$ at the reader's antenna is assumed to have zero mean with $\mathbb{E}\left[|n_r|^2\right] = \sigma^2$. Further, the modulation index $m_{i,k} \in \left[0,1\right]$ is the characteristic difference between the two transmit symbols $S_i$ and $S_k$,and is defined below~\cite{fuschini2008analytical}:
\begin{align}
    m_{i,k} &\triangleq \frac{\abs{\Gamma_i - \Gamma_k}}{2} \nonumber \\
      &\stackrel{(s_3)}{=} \frac{\sqrt{\left(\Gamma_{{\rm a}i}-\Gamma_{{\rm a}k}\right)^2 + \left(\Gamma_{{\rm b}i}-\Gamma_{{\rm b}k}\right)^2}}{2},  
      \label{eq6}
\end{align}
$\forall i,k \in \mathbb{M}, \hspace{1mm} i \neq k,$ and $(s_3)$ is obtained using \eqref{eq3}. We set $\nu = \frac{\abs{V_0} m_{i,k}}{2 \sqrt{2} \sigma}$, and the complementary error function $\text{erfc}\left( \nu \right) = 1 - \text{erf}\left( \nu \right)$. Given that $\text{erf} \left( \nu \right)$ is the error function, this implies the higher the $m_{i,k}$, the lower the SER. In other words, each transmitted symbol must exhibit unique characteristics to be distinguished from other symbols. This can be achieved by ensuring that the modulation index between any two transmitted symbols exceeds a threshold $m_{\rm th}$. Thus, it is crucial to set a design requirement of $m_{i,k} \geq m_{\rm th}$ to keep the SER within acceptable limits.

\section{Problem Definition} \label{sectionPD}

 This section presents the problem formulation for maximizing the average harvested power at the passive tag in two distinct cases. In Section~\ref{IV.A}, we consider the general case where the BackCom system operates with $M$-ASK modulation scheme and satisfies the underlying reader sensitivity, tag sensitivity, and SER requirements. We formulate an optimization problem to determine the optimal reflection coefficients that maximize the average harvested power. In the second case, we formulate a new problem by considering the BASK modulation scheme and relaxing the reader sensitivity constraint. Specifically, we assume that the $P_{{\rm b}i} \geq P_{\rm b,min}$ is always satisfied. Further details regarding this assumption are provided in Section~\ref{IV.B}.

\subsection{Optimization Formulation for the General Case}\label{IV.A}
 
 When activated, the tag arbitrarily backscatters the symbol $S_i$ that carries $\log_2 M$ number of bits while delivering $P_{{\rm L}i}$ to the tag. We denote $p_i$ as the occurrence probability of $S_i$, with $p_i \in \left[0,1\right]$ and $\sum_{i=1}^M p_i = 1$. Furthermore, we set $\mathbf{p} = \left[ p_1, p_2, ..., p_M \right]$, and assume $p_i \geq p_{i+1}$ without any loss of generality. In general, it is not necessary for all transmit symbols to have equal probabilities of occurrence. Non-equiprobable signaling has been explored and can be advantageous~\cite{yang2013non,wei2012optimized}. Hence, we consider $p_i$ as application dependent constant, and the average harvested power $P_{\rm L,avg}$ at the tag is given by:
\begin{equation}
     P_{\rm L,avg} = \smashoperator[r]{\sum_{i=1}^M} p_i P_{{\rm L}i} ,
     \label{eq7}
 \end{equation}

 Maximizing $P_{\rm L,avg}$ is motivated by the need to provide more power to the tag, which is crucial for low-power and low-bandwidth applications such as backscatter-assisted IoT sensing and identification.\footnote{For a tag operating with HTT protocol, the total energy consumed for circuitry operation during backscattering includes both the energy harvested during the energy harvesting phase and the backscattering phase. As a result, a greater $P_{\rm L,avg}$ shortens the duty cycle, allowing for increased data transmission.} The increase in $P_{\rm L,avg}$ at the tag allows it to power additional sensing units, such as temperature sensors, humidity sensors, accelerometers, and cameras. This, in turn, enables the collection of more data using the same resources for BackCom. Hence, maximizing $P_{\rm L,avg}$ is essential for optimizing the efficiency of the IoT network.

 Since $P_{\rm L,avg}$ depends on $\Gamma_i$, the goal is to determine the optimal reflection coefficients that maximize $P_{\rm L,avg}$, subject to the following constraints. Constraint $C1$ defines the domain of the reflection coefficient $\abs{\Gamma_i} \leq 1$, while $C2$ and $C3$ establish the boundary conditions for $\Gamma_{{\rm a}i}$ and $\Gamma_{{\rm b}i}$, respectively. To meet the minimum SER requirement, the tag's modulation index $m_{i,k}$ must exceed $m_{\rm th}$, as specified in constraint $C4$. Additionally, constraint $C5$ ensures continuous operation by requiring that $P_{{\rm L}i} \geq P_{\rm L,min}$. Finally, constraint $C6$ requires the backscattered signal of all different symbols must exceed the minimum backscattered power threshold $P_{\rm b,min}$, ensuring that the reader can accurately receive and decode the signal. With these constraints in place and maximizing $P_{\rm L,avg}$ as the objective, the optimization problem is defined as:
\begin{align*}
(P1): &\max_{\mathbf{\Gamma}} \quad P_{\rm L,avg} \\
\text{s.t.} \hspace{1.5mm}
    C1:& \Gamma^2_{{\rm a}i} + \Gamma^2_{{\rm b}i} \leq 1, \hspace{2mm} \forall i \in \mathbb{M},\\
    C2:& \Gamma_{{\rm a}i} \in  \left[-1,1\right], \hspace{2mm} \forall i \in \mathbb{M},\\
    C3:& \Gamma_{{\rm b}i} \in  \left[-1,1\right], \hspace{2mm} \forall i \in \mathbb{M},\\
    C4:& \frac{\sqrt{\left(\Gamma_{{\rm a}i}-\Gamma_{{\rm a}k}\right)^2 + \left(\Gamma_{{\rm b}i}-\Gamma_{{\rm b}k}\right)^2}}{2} \geq m_{\rm th}, \\ 
    & \hspace{5 mm} \forall i,k \in \mathbb{M}, \hspace{1mm} i \neq k, \\
    C5:& \eta_{\rm h} P_{\rm a} \left(1-\Gamma_{{\rm a}i}^2-\Gamma_{{\rm b}i}^2\right) \geq \Delta P_{\rm L,min}, \forall i \in \mathbb{M},\\
    C6:& \eta_{\rm b} P_{\rm a} G_{\rm t} \Bigl[ \left( 1 - \Gamma_{{\rm a}i} \right)^2 + \Gamma^2_{{\rm b}i} \Bigr] \geq P_{\rm b,min}, \hspace{2mm} \forall i \in \mathbb{M},
\end{align*}
where $\mathbf{\Gamma} = \left[\Gamma_{\rm a1},\Gamma_{\rm b1},\Gamma_{\rm a2},\Gamma_{\rm b2},...,\Gamma_{{\rm a}M},\Gamma_{{\rm b}M}\right]$ are the optimization variables and $\Delta P_{\rm L,min} = P_{\rm L,min} + \eta_{\rm h} P_{\rm sen}$. By solving $\left(P1\right)$, we determine the maximum $P_{\rm L,avg}$ $\left(\right.$denoted as $P^\ast_{\rm L,avg}\left.\right)$ by jointly optimizing $\mathbf{\Gamma}$ (denoted as $\mathbf{\Gamma}^\ast$). To simplify the problem-solving process, we introduce the following lemmas to reduce problem $(P1)$ from a $2M$-variable optimization problem to an $M$-variable problem.


 \begin{lemma}\label{L1}
 The average harvested power is maximized when either $\Gamma_{{\rm a}i}=0$ or $\Gamma_{{\rm b}i}=0$, $\forall i \in \mathbb{M}$.
 \end{lemma}

 \begin{proof}
 First, on setting $\Gamma_{{\rm b}i}$ as constant, we find that $\frac{\partial^2 P_{{\rm L}i}}{\partial \Gamma^2_{{\rm a}i}} = - 2 \eta_{\rm h} P_{\rm a}$, which implies $P_{{\rm L}i}$ is a concave function in $\Gamma_{{\rm a}i}$. Therefore, for a given $\Gamma_{{\rm b}i}$, the optimal value of $\Gamma_{{\rm a}i}$ maximizing $P_{{\rm L}i}$ as obtained by solving $\frac{\partial P_{{\rm L}i}}{\partial \Gamma_{{\rm a}i}} =0$ is $\Gamma_{{\rm a}i}=0$. Likewise, we set $\Gamma_{{\rm a}i}$ as a constant and we find that $\frac{\partial^2 P_{{\rm L}i}}{\partial \Gamma^2_{{\rm b}i}} = - 2 \eta_{\rm h} P_{\rm a}$, which implies $P_{{\rm L}i}$ is also a concave function in $\Gamma_{{\rm b}i}$. Similarly, for a given $\Gamma_{{\rm a}i}$, the maximum harvested power as obtained by solving $\frac{\partial P_{{\rm L}i}}{\partial \Gamma_{{\rm b}i}} =0$ is $\Gamma_{{\rm b}i}=0$. Hence, we proved Lemma 1.
\end{proof}

\begin{lemma}\label{L2}
Given the same harvested power, a real reflection coefficient, where $\Im(\Gamma_i) = \Gamma_{{\rm b}i} = 0$, ensures better receiver sensitivity at the reader compared to a purely imaginary reflection coefficient, where $\Re(\Gamma_i) = \Gamma_{{\rm a}i} = 0$.
\end{lemma} 

 \begin{proof}
 Refer to Appendix \ref{ApdA} for the proof of Lemma 2. 
 \end{proof}

Using Lemmas \ref{L1} and \ref{L2}, we can reformulate problem $\left(P1\right)$ into the following problem $\left(P2\right)$: 
\begin{align*}
&(P2): \max_{\mathbf{\Gamma_a}} \quad P_{\rm L,avg} \\
\text{s.t.} \hspace{3mm} &C2, C7: \frac{\abs{\Gamma_{{\rm a}i}-\Gamma_{{\rm a}k}}}{2} \geq m_{\rm th},  \hspace{2mm} \forall i,k \in \mathbb{M}, \hspace{1mm} i \neq k,\\
        &C8: \eta_{\rm h} P_{\rm a} \left(1-\Gamma_{{\rm a}i}^2\right) \geq \Delta P_{\rm L,min}, \hspace{2mm} \forall i \in \mathbb{M},\\
        &C9: \eta_{\rm b} P_{\rm a} G_{\rm t} \left( 1 - \Gamma_{{\rm a}i} \right)^2 \geq P_{\rm b,min}, \hspace{2mm} \forall i \in \mathbb{M},
\end{align*}
where $\mathbf{\Gamma_a} = \left[ \Gamma_{\rm a1}, \Gamma_{\rm a2},..., \Gamma_{{\rm a}M} \right]$, and $\left(P2\right)$ is a $M-$variable problem. By solving problem $\left(P2\right)$, we will obtain $P^\ast_{\rm L,avg}$ along with the optimal $\mathbf{\Gamma_a}$, denoted as $\mathbf{\Gamma^\ast_{\rm a}} = \left[ \Gamma^\ast_{\rm a1}, \Gamma^\ast_{\rm a2}, \ldots, \Gamma^\ast_{{\rm a}M} \right]$.

\subsection{Optimization with BASK Modulation}\label{IV.B}

Here, we study the binary encoding scheme used in backscatter tags, known as BASK modulation, which is widely used in RFID applications. BASK can also be referred to as 2-ASK, where $M=2$. The transmit symbols contain only one bit, either ‘1’ or ‘0’. As a result, the probabilities for these 2 symbols are $p_1$ and $p_2$, with $p_1 + p_2 = 1$. The average harvested power at the tag with BASK modulation is expressed as: 
\begin{equation}
    \hat{P}_{\rm L,avg} = p_1 P_{\rm L1} + \left(1-p_1\right) P_{\rm L2}
\end{equation}

Likewise, we aim to maximize $\hat{P}_{\rm L,avg}$ by optimizing $\hat{\mathbf{\Gamma}} = \left[ \Gamma_{\rm a1}, \Gamma_{\rm b1}, \Gamma_{\rm a2}, \Gamma_{\rm b2} \right]$ under the same BackCom system operational requirements outlined in problem $\left(P1\right)$. In contrast, we only consider $C1-C5$ and assume $C6$ is always satisfied in BackCom systems. This assumption is valid because commercial readers typically exhibit good receive sensitivity, and RFID tags are typically downlink-limited.\footnote{The commercial Impinj RFID reader has a receive sensitivity of $-84$ dBm~\cite{impinj2023}. Additionally, in monostatic BackCom systems, the distance between the reader and the tag is generally short, as the passive tag relies on harvesting power from the reader's transmitted RF carrier in the downlink. As a result, the minimum backscattered power constraint is less critical and can often be relaxed compared to the more important tag sensitivity constraint.} The corresponding optimization problem is then formulated as follows:
\begin{align*}
(P3): &\max_{\hat{\mathbf{\Gamma}}} \quad \hat{P}_{\rm L,avg}  \\
\text{s.t.} \hspace{3mm} &C1, C2, C3, C4, C5.
\end{align*}

Similarly, we apply Lemmas \ref{L1} and \ref{L2} to reduce the 4-variable problem $\left(P3\right)$ to a 2-variable problem by setting $\Gamma_{\rm b1}=\Gamma_{\rm b2}=0$. Without loss of generality, we assume $\Gamma_{\rm a1} \geq \Gamma_{\rm a2}$, and problem $\left(P3\right)$ is reformulated as:
 \begin{align*}
&(P4): \max_{\Gamma_{\rm a1},\Gamma_{\rm a2}} \quad \hat{P}_{\rm L,avg} \\
\text{s.t.} \hspace{3mm} &C2, C10: \frac{\Gamma_{\rm a1}-\Gamma_{\rm a2}}{2} \geq m_{\rm th},\\
        &C11: \eta_{\rm h} P_{\rm a} \left(1-\Gamma_{{\rm a}i}^2\right) \geq \Delta P_{\rm L,min}, \hspace{2mm} \forall i \in \{1,2\}.
\end{align*}

 To distinguish the results of different formulated problems, we denote the maximum average harvested power obtained from problem $(P4)$ as $\hat{P}^\ast_{\rm L,avg}$, determined using the optimal solution $\hat{\mathbf{\Gamma}}^\ast = \left[ \hat{\Gamma}^\ast_{\rm a1}, \hat{\Gamma}^\ast_{\rm a2} \right]$. The methods for solving problems $(P2)$ and $(P4)$ are detailed in Section~\ref{Solution}.
 
 \section{Solution Methodology}\label{Solution}
 \subsection{Problem Feasibility}\label{Sec:PFaC}

Before solving problems $\left(P2\right)$ and $\left(P4\right)$, we will need to verify their feasibility with the given $\Delta P_{\rm L,min}$, $P_{\rm b,min}$ and $m_{\rm th}$. This step is crucial to ensure that an optimal set of reflection coefficients exists and meets all the operational requirements of the BackCom system. The following feasibility check applies to problem $\left(P2\right)$ and can be extended to problem $\left(P4\right)$ if we relax $C9$ as it does not have $P_{\rm b,min}$ constraint. 

First, considering the unconstrained problem $\left(P2\right)$, we notice that the $P_{\rm L,avg}$ is maximized when $\Gamma_{{\rm a}i} = 0$, $\forall i \in \mathbb{M}$. However, constraint $C7$ does not allow this to happen because it sets the separation between any two reflection coefficients $\abs{\Gamma_{{\rm a}i}-\Gamma_{{\rm a}k}}$ must be greater than $2 m_{\rm th}$. Since the separation between any two reflection coefficients is different, the minimum separation is specified in the following Lemma. 

\begin{lemma}\label{L123}
    The minimum separation between the reflection coefficients from the optimal solution set is $\abs{\Gamma^\ast_{{\rm a}i}-\Gamma^\ast_{{\rm a}k}} = 2 m_{\rm th}$, which is determined by setting them to satisfy $C7$ at equality.
\end{lemma}

\begin{proof}
    We first set a reflection coefficient at the perfect matched condition, denoted as $\grave{\Gamma}^{\left(0\right)} = 0$. Then, we set 2 different real reflection coefficients $\grave{\Gamma}^{\left(1\right)}$ and $\grave{\Gamma}^{\left(2\right)}$, and the modulation indices between these 2 reflection coefficients and $\grave{\Gamma}^{\left(0\right)}$ are $m_1 = \frac{\abs{\grave{\Gamma}^{\left(1\right)} - \grave{\Gamma}^{\left(0\right)}}}{2}$ and $m_2 = \frac{\abs{\grave{\Gamma}^{\left(2\right)} - \grave{\Gamma}^{\left(0\right)}}}{2}$, respectively. As we set $m_1 < m_2$, we obtain: 
    \begin{equation}
        \abs{\grave{\Gamma}^{\left(1\right)}} < \abs{\grave{\Gamma}^{\left(2\right)}}. 
        \label{eq456}
    \end{equation}
    Since we know, from \eqref{eq4}, the harvested power is greater with a lower reflection coefficient magnitude under the same operating circumstance. Hence, $\grave{\Gamma}^{\left(1\right)}$ will give a greater harvested power as compare with $\grave{\Gamma}^{\left(2\right)}$. This indicates that the minimum separation between reflection coefficients must be the lowest to achieve greater $P_{\rm L,avg}$. Since constraint $C7$ needs to be satisfied, the minimum separation between the optimal reflection coefficients will satisfy it at equality, such that $\abs{\Gamma^\ast_{{\rm a}i}-\Gamma^\ast_{{\rm a}k}} = 2 m_{\rm th}$. Hence, we proved Lemma~\ref{L123}. 
\end{proof}

 Lemma \ref{L123} yields a corollary as stated below: 

\begin{corollary}
 The difference between the largest optimal reflection coefficient $\Gamma_{L}$ and smallest optimal reflection coefficient $\Gamma_{S}$, in $\mathbf{\Gamma^\ast_a}$ is given by $\Gamma_{L} - \Gamma_{S} = 2 m_{\rm th} \left(M-1\right)$.
\end{corollary}

\begin{proof}
    Using Lemma \ref{L123} and knowing $\Gamma_{{\rm a}i}$ is a real number, we can represent all $\Gamma_{{\rm a}i}$, $\forall i = \mathbb{M}$, on a number line with an equal gap of $2 m_{\rm th}$ for any two adjacent reflection coefficients. Since there is $M$ number of reflection coefficients in the $M$-ASK modulated tag design, there are $\left(M-1\right)$ gaps between the largest and smallest optimal reflection coefficients. Hence, $\Gamma_{L} - \Gamma_{S} = 2 m_{\rm th} \left(M-1\right)$. 
\end{proof}

From constraints $C8$ and $C9$, we can determine the upper bound $\Gamma_{\rm ub}$ and lower bound $\Gamma_{\rm lb}$ of $\Gamma_{{\rm a}i}$. By rearranging $C8$ and $C9$, we able to show $\Gamma_{{\rm a}i}$ lies in the following range, respectively: 
\begin{equation}
    -\sqrt{1-\frac{\Delta P_{\rm L,min}}{ \eta_{\rm h} P_{\rm a}}} \leq \Gamma_{{\rm a},i} \leq \sqrt{1-\frac{\Delta P_{\rm L,min}}{ \eta_{\rm h} P_{\rm a}}},
    \label{eqC8}
\end{equation}

\begin{equation}
    \Gamma_{{\rm a},i} \leq 1 - \sqrt{\frac{P_{\rm b,min}}{ \eta_{\rm b} P_{\rm a} G_{\rm t}}} , \hspace{3mm} \Gamma_{{\rm a},i} \geq 1 + \sqrt{\frac{P_{\rm b,min}}{\eta_{\rm b} P_{\rm a} G_{\rm t}}}.
    \label{eqC9}
\end{equation}

From \eqref{eqC8} and \eqref{eqC9} along with constraint $C2$, $\Gamma_{\rm ub}$ and $\Gamma_{\rm lb}$ are determined as: 
\begin{align}
    \Gamma_{\rm ub} &= \min \Biggl\{ {\sqrt{1-\frac{\Delta P_{\rm L,min}}{ \eta_{\rm h} P_{\rm a}}}, 1 - \sqrt{\frac{P_{\rm b,min}}{ \eta_{\rm b} P_{\rm a} G_{\rm t}}}} \hspace{1mm} \Biggr\} \label{eq13b}, \\
    \Gamma_{\rm lb} &= -\sqrt{1-\frac{\Delta P_{\rm L,min}}{ \eta_{\rm h} P_{\rm a}}} \label{eq14}.
\end{align}

Subsequently, we form a condition with the defined $\Gamma_{\rm ub}$ and $\Gamma_{\rm lb}$ to verify problem $\left(P2\right)$ has at least one possible solution set. The feasibility check condition is specified in Theorem \ref{T1}.

\begin{theorem}\label{T1}
    With the selected order of $M$-ASK modulation scheme, the formulated problems are feasible if and only if $M \leq 1 + \frac{\Gamma_{\rm ub} - \Gamma_{\rm lb}}{2 m_{\rm th}}$.
\end{theorem}

\begin{proof}
    Since we know $\Gamma_{L} - \Gamma_{S}$ is the range that covers all $\Gamma_{{\rm a}i}$ on a number line, and we have proved that the upper and lower boundaries for $\Gamma_{{\rm a}i}$ are respectively $\Gamma_{\rm ub}$ and $\Gamma_{\rm lb}$, $\Gamma_{L} - \Gamma_{S}$ must fit within the boundary to ensure feasibility. Thus, we have the condition $\Gamma_{L} - \Gamma_{S} \leq \Gamma_{\rm ub} - \Gamma_{\rm lb}$. By substituting $\Gamma_{L} - \Gamma_{S} = 2 m_{\rm th} \left(M-1\right)$ and rearranging, we will obtain $M \leq 1 + \frac{\Gamma_{\rm ub} - \Gamma_{\rm lb}}{2 m_{\rm th}}$.
\end{proof}

\begin{remark}
    From Theorem~\ref{T1}, we know that $m_{\rm th} \leq \frac{\Gamma_{\rm ub} - \Gamma_{\rm lb}}{2\left(M-1\right)}$, indicating that the highest value of $m_{\rm th}$ decreases with increasing $M$. Consequently, the minimum achievable SER decreases with $M$.
\end{remark}
\vspace{-4mm}

 \subsection{Optimal Reflection Coefficients for BASK Modulation}\label{SecB}

We first solve the relatively simpler problem $\left(P4\right)$ to get a basic idea of the optimal reflection coefficients selection before solving $\left(P2\right)$. We start with invoking a key characteristic of problem $\left(P4\right)$ using Lemma \ref{L4}.

 \begin{lemma}\label{L4}
 Problem $\left(P4\right)$ is a convex problem.
 \end{lemma}
 
 \begin{proof}
 Refer to Appendix \ref{ApdB} for the proof of Lemma \ref{L4}.
 \end{proof}

\vspace{2mm}
 
 Since $\left(P4\right)$ is a convex problem, we can claim that the KKT point is the global optimal solution. Keeping constraint $C2$ implicit, the Lagrangian of $\left(P4\right)$ can be written as:
 \begin{align}
     \mathcal{L}_1 &= -p_1 \eta_{\rm h} P_{\rm a} \left(1-\Gamma_{\rm a1}^2\right) - \left(1-p_1\right) \eta_{\rm h} P_{\rm a} \left(1-\Gamma_{\rm a2}^2\right) \nonumber \\ &+ \lambda_1 \left( m_{\rm th} - \frac{\Gamma_{\rm a1}-\Gamma_{\rm a2}}{2} \right) + \lambda_2 \left(\Gamma^2_{\rm a1} - 1 + \frac{\Delta P_{\rm L,min}}{ \eta_{\rm h} P_{\rm a}}\right) \nonumber \\
     &+ \lambda_3 \left(\Gamma^2_{\rm a2} - 1 + \frac{\Delta P_{\rm L,min}}{ \eta_{\rm h} P_{\rm a}}\right),
     \label{eq8}
 \end{align}
where $\lambda_1$ represents the Lagrange multiplier associated with $C10$, and $\lambda_2,\lambda_3$ correspond to $C11$. The KKT point can be found by solving the following equations.
\vspace{5mm}
\begin{align}
    \frac{\partial \mathcal{L}_1}{\partial \Gamma_{\rm a1}} = 2 p_1 \eta_{\rm h} P_{\rm a} \Gamma_{\rm a1} - \frac{1}{2} \lambda_1 + 2 \lambda_2 \Gamma_{\rm a1} &= 0\label{eq9},\\
    \frac{\partial \mathcal{L}_1}{\partial \Gamma_{\rm a2}} = 2 \left(1-p_1\right) \eta_{\rm h} P_{\rm a} \Gamma_{\rm a2} + \frac{1}{2} \lambda_1 + 2 \lambda_3 \Gamma_{\rm a2} &= 0,\label{eq10}\\
    \lambda_1 \left( m_{\rm th} - \frac{\Gamma_{\rm a1}-\Gamma_{\rm a2}}{2} \right) &= 0,\label{eq11}\\
    \lambda_2 \left(\Gamma^2_{\rm a1} - 1 + \frac{\Delta P_{\rm L,min}}{ \eta_{\rm h} P_{\rm a}} \right) &= 0,\label{eq12}\\
    \lambda_3 \left(\Gamma^2_{\rm a2} - 1 + \frac{\Delta P_{\rm L,min}}{ \eta_{\rm h} P_{\rm a}} \right) &= 0.\label{eq13}
\end{align}
along with $\lambda_1, \lambda_2, \lambda_3 \geq 0$. Equations \eqref{eq9} and \eqref{eq10} are the sub-gradient conditions, and \eqref{eq11},\eqref{eq12},\eqref{eq13} are the complementary slackness conditions. While solving \eqref{eq9} $-$ \eqref{eq13}, we obtain $\mathbf{\hat{\Gamma}^\ast}$ in terms of the constant parameters, and thereby determine $\hat{P}^\ast_{\rm L,avg}$. Next, we discuss the method to determine the underlying KKT point $\left( \mathbf{\hat{\Gamma}^\ast}, \lambda^\ast_1, \lambda^\ast_2, \lambda^\ast_3 \right)$ in Lemma \ref{L5}.

\begin{lemma}\label{L5}
The $\mathbf{\hat{\Gamma}^\ast}$ is obtained from one of the 3 cases, which are case $\left(a\right): \lambda_1 \neq 0, \lambda_2=\lambda_3=0$, case $\left(b\right): \lambda_1 \neq 0, \lambda_2 \neq 0, \lambda_3=0$, and case $\left(c\right): \lambda_1 \neq 0, \lambda_3 \neq 0, \lambda_2=0$. 
\end{lemma}

\begin{proof}
First, we check whether $\lambda_1$ can be zero or not. Since the objective function is decreasing with $\Gamma_{\rm a1}$ and $\Gamma_{\rm a2}$, the optimal solution without the constraints will have $\Gamma_{\rm a1}=\Gamma_{\rm a2}=0$. However, constraint $C10$ requires a minimum separation between $\Gamma_{\rm a1}$ and $\Gamma_{\rm a2}$. Therefore, as proved, constraint $C10$ is satisfied at equality, which implies $\lambda_1$ is always positive. 

It is noticed that $\lambda_2$ and $\lambda_3$ simultaneously greater than zero is very unlikely because this is feasible only when $m_{\rm th} = \sqrt{1 - \frac{\Delta P_{\rm L,min}}{\eta_{\rm h} P_{\rm a}}}$. When this condition holds, $\mathbf{\hat{\Gamma}^\ast}$ can be obtained from either case $\left(b\right)$ or case $\left(c\right)$. Therefore, the optimal solution is given by either or both $\lambda_2$ and $\lambda_3$ are zero, while $\lambda_1 > 0$.
\end{proof}

Thus, 3 different $\Gamma_{\rm a1}$ can be obtained from the 3 cases stated in Lemma \ref{L5}, in which we denote $\Gamma_{\rm a1}$ obtained from case $\left( a \right)$, $\left( b \right)$ and $\left( c \right)$ as $\hat{\Gamma}^{\left(\rm a\right)}_{\rm a1}$, $\hat{\Gamma}^{\left(\rm b\right)}_{\rm a1}$, and $\hat{\Gamma}^{\left(\rm c\right)}_{\rm a1}$, respectively. In case $\left( a \right)$ where $\lambda_2=\lambda_3=0$, $\hat{\Gamma}^{\left(\rm a\right)}_{\rm a1}$ is obtained with \eqref{eq9}, \eqref{eq10} and \eqref{eq11}. On the other hand, in case $\left( b \right)$ where $\lambda_3=0$, $\hat{\Gamma}^{\left(\rm b\right)}_{\rm a1}$ is obtained with \eqref{eq12}. Likewise, we use \eqref{eq11} and \eqref{eq13} to obtain $\hat{\Gamma}^{\left(\rm c\right)}_{\rm a1}$ in case $\left( c \right)$ where $\lambda_2=0$. As a result, the closed-form expression of the solutions are:
\begin{subequations}
\begin{align}
    \hat{\Gamma}^{\left(\rm a\right)}_{\rm a1} &= 2 \left(1-p_1\right) m_{\rm th}, \\
    \hat{\Gamma}^{\left(\rm b\right)}_{\rm a1} &= \sqrt{1 - \frac{\Delta P_{\rm L,min}}{ \eta_{\rm h} P_{\rm a}}}, \\
    \hat{\Gamma}^{\left(\rm c\right)}_{\rm a1} &= -\sqrt{1 - \frac{\Delta P_{\rm L,min}}{ \eta_{\rm h} P_{\rm a}}}+2 m_{\rm th}.
\end{align}
\label{eqE22}
\end{subequations}

From Lemma \ref{L4}, we recall that problem $\left(P4\right)$ is convex. Therefore, only one KKT point will satisfy all the constraints, which implies $\hat{\Gamma}^\ast_{\rm a1}$ is selected from $\hat{\Gamma}^{\left(\rm a\right)}_{\rm a1}$, $\hat{\Gamma}^{\left(\rm b\right)}_{\rm a1}$ and $\hat{\Gamma}^{\left(\rm c\right)}_{\rm a1}$. Notably, all the optimal Lagrange multipliers $\left( \lambda^\ast_1, \lambda^\ast_2, \lambda^\ast_3 \right)$ of the KKT point must be positive for it to be the optimal solution.  It is noticed that $\lambda_2$ from case $\left(b\right)$ is positive only when $\hat{\Gamma}^{\left(\rm b\right)}_{\rm a1} < \hat{\Gamma}^{\left(\rm a\right)}_{\rm a1} $, whereas $\lambda_3$ from case $\left(c\right)$ is positive only when $\hat{\Gamma}^{\left(\rm c\right)}_{\rm a1} > \hat{\Gamma}^{\left(\rm a\right)}_{\rm a1}$. Besides, from \eqref{eq13b}, \eqref{eq14} and Lemma \ref{L123}, we notice that the range of $\hat{\Gamma}^\ast_{\rm a1}$ is given as $\left[ \hat{\Gamma}^{\left(\rm c\right)}_{\rm a1}, \hat{\Gamma}^{\left(\rm b\right)}_{\rm a1} \right]$. Hence, incorporating all these properties, $\hat{\Gamma}^\ast_{\rm a1}$ is determined by:
\begin{equation}
    \hat{\Gamma}^\ast_{\rm a1} = \max \biggl\{ \hat{\Gamma}^{\left(\rm c\right)}_{\rm a1}, \min \Bigl\{ \hat{\Gamma}^{\left(\rm a\right)}_{\rm a1}, \hat{\Gamma}^{\left(\rm b\right)}_{\rm a1} \Bigr\} \biggr\}, 
    \label{eqE23}
\end{equation}
and $\hat{\Gamma}^\ast_{\rm a2}$ is then determined by:
\begin{equation}
    \hat{\Gamma}^\ast_{\rm a2} = \hat{\Gamma}^\ast_{\rm a1} - 2 m_{\rm th}.
\end{equation}

Eventually, the maximum average harvested power that is determined from problem $\left(P4\right)$ is: 
\begin{align}
    \hat{P}^\ast_{\rm L,avg} &= \eta_{\rm h} P_{\rm a} \Bigg[ p_1 \left(1-\left(\hat{\Gamma}^{\ast}_{\rm a1}\right)^2\right) + p_2 \left(1-\left(\hat{\Gamma}^{\ast}_{\rm a2}\right)^2\right) \Bigg]. \nonumber \\
    & 
    \label{eqE25}
\end{align}

\vspace{-3mm}
Remarkably, we obtain the closed-form expression for the global solution of problem $\left(P4\right)$. This signifies the utility of using KKT conditions to solve the problem. From~\eqref{eqE22}, \eqref{eqE23}, and~\eqref{eqE25}, it is evident that the occurrence probabilities $p_1$ and $p_2$ are crucial factors. The maximum average harvested power in the unequal symbol probability scenario is greater than in the equal symbol probability scenario, as illustrated graphically in Sec.~\ref{sectionRD}.

\begin{remark}
    In the context of this case, it is observed that $\hat{\Gamma}^\ast_{{\rm a}i}$ is closely linked to $p_i$. In particular, $\hat{\Gamma}^\ast_{{\rm a}i}$ will be set closer to zero when $p_i$ is higher. This highlights the significance of $p_i$ in determining $\hat{\Gamma}^\ast_{{\rm a}i}$. The implications of this observation can contribute to a better understanding of the underlying system and inform the development of more efficient and effective strategies for optimization.
\end{remark}

\subsection{Optimal Reflection Coefficients for General Case}\label{secSolution}

\begin{figure}[!t]
    \centering
    \includegraphics[width=3.3in]{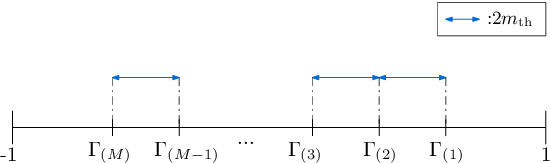}
    \caption{The descending order sequence of $\Gamma_{\left(k\right)}$ on number line}
    \label{Fig12345}
\end{figure}

 Here, we propose a method to solve problem $\left(P2\right)$ and obtain the optimal reflection coefficients for the tag with $M$-ASK modulation in the general case. It is evident that $\left(P2\right)$ is a non-convex problem as constraints $C7$ and $C9$ are non-convex. Therefore, we cannot obtain the global solution by applying the KKT conditions unless we transform $\left(P2\right)$ into a convex problem. Next, we will show how the original complex problem $\left(P2\right)$ can be solved in 2 steps. Specifically, we first transform the original problem $\left(P2\right)$ into a new convex problem under an assumption that will allow us to obtain a candidate of the global solution with the KKT conditions. Later, we determine all the candidates where the optimal solution is the solution set that gives the maximum $P_{\rm L,avg}$. The details of these 2 steps are discussed in the following paragraphs.  
 
 Following Lemma \ref{L123}, we can sequence elements of $\mathbf{\Gamma_a}$ on a number line, with the gap between any 2 adjacent reflection coefficients being $2 m_{\rm th}$. Although we still do not know the sequence of all $\Gamma^\ast_{{\rm a}i}$, this characteristic will hold when we locate them on a number line. In this context, we denote $\Gamma_{\left(k\right)}$, $\forall k \in \mathbb{M}$ to represent the placement of all $\Gamma_{{\rm a}i}$ in descending order as shown in Fig.~\ref{Fig12345}. Further, we set $p_{\left(k\right)}$ as the probability corresponding to $\Gamma_{\left(k\right)}$, and $\mathbf{p}_2 = \left[ p_{\left(1\right)}, p_{\left(2\right)},..., p_{\left(M\right)} \right]$. This implies the sequence of all $\Gamma_{{\rm a}i}$ can be determined with the corresponding $p_i$ if all $p_{\left(k\right)}$ are known. 

Assuming the sequence of all $\Gamma^\ast_{{\rm a}i}$ is known, which indicates $p_{\left(k\right)}$ is also known, we can use $\Gamma_{\left(k\right)}$, $\forall k \in \mathbb{M}$ as the optimization variables to determine $P^\ast_{\rm L,avg}$. This assumption will be further validated by a proposed approach discussed later in the same section. With this assumption, we transform problem $\left(P2\right)$ into problem $\left(P5\right)$, where constraint $C7$ can be equivalently replaced by $C12$, which ensures the gap between any 2 adjacent $\Gamma_{\left(k\right)}$ is $2 m_{\rm th}$. Besides, as stated earlier, constraints $C2$, $C8$, and $C9$ provide the upper bound $\Gamma_{\rm ub}$ and lower bound $\Gamma_{\rm lb}$ for $\Gamma_{{\rm a}i}$, which we can combine all these 3 constraints and form a new constraint $C13$ as $\Gamma_{\rm lb} 
\leq \Gamma_{\left(k\right)} \leq \Gamma_{\rm ub}$. Hence, under the assumption of $p_{\left(k\right)}$ is known, problem $\left(P2\right)$ is equivalent to:
\begin{align*}
(P5): &\max_{\mathbf{\Gamma_{\left(\rm a\right)}}} \quad P_{\left(\rm L,avg\right)} = \smashoperator[r]{\sum^M_{k=1}} p_{\left(k\right)} \eta_{\rm h} P_{\rm a} \Bigl( 1- \Gamma^2_{\left(k\right)} \Bigr) \\
\text{s.t.} \hspace{3mm} &C12: \Gamma_{\left(k\right)} = \Gamma_{\left(1\right)} - 2 m_{\rm th} \left(k-1\right), \hspace{2mm} \forall k \in \mathbb{M}, \\
    &C13: \Gamma_{\rm lb} \leq \Gamma_{\left(k\right)} \leq \Gamma_{\rm ub}, \hspace{2mm} \forall k \in \mathbb{M},
\end{align*}
where $\mathbf{\Gamma_{\left(\rm a\right)}} = \left[ \Gamma_{\left(1\right)},\Gamma_{\left(2\right)},...,\Gamma_{\left(M\right)} \right]$. We then substitute $\Gamma_{\left(k\right)}$ from constraint $C12$ into the objective function and turn it into a single-variable function of $\Gamma_{\left(1\right)}$, as follows:
\begin{align}
      &P_{\left(\rm L,avg\right)} = \eta_{\rm h} P_{\rm a} \Biggl[ 1- \smashoperator[r]{\sum^M_{k=1}} p_{\left(k\right)} \Bigl( \Gamma_{\left(1\right)} - 2 m_{\rm th} \left(k-1\right) \Bigr)^2 \Biggr] \nonumber \\
      &
\end{align}

\vspace{-4mm}
Additionally, since we know $\Gamma_{\left(k\right)} > \Gamma_{\left(k+1\right)}$, constraint $C13$ will always satisfy when we set $\Gamma_{\left(1\right)} \leq \Gamma_{\rm ub}$ and $\Gamma_{\left(M\right)} \geq \Gamma_{\rm lb}$. Consequently, problem $(P5)$ can be further transformed into a single variable problem as follows: 
\begin{align*}
(P6): &\max_{\Gamma_{\left(1\right)}} \quad P_{\left(\rm L,avg\right)} \\
\text{s.t.} \hspace{3mm} &C14: \Gamma_{\left(1\right)} - 2 m_{\rm th} \left(M-1\right) \geq \Gamma_{\rm lb}, \\
    &C15: \Gamma_{\left(1\right)} \leq \Gamma_{\rm ub},
\end{align*}
where $C14$ and $C15$ are the boundary conditions of $\Gamma_{\left(1\right)}$. Next, Lemma~\ref{L66} discusses the convexity of problem $\left(P6\right)$.

\begin{lemma}\label{L66}
    Problem $\left(P6\right)$ is a convex problem.
\end{lemma}

\begin{proof}
    Firstly, we determine the second derivative of the objective function, which is $\frac{\partial^2 P_{\left(\rm L,avg\right)}}{\partial \Gamma_{\left(1\right)}^2} = - 2 \eta_{\rm h} P_{\rm a} \smashoperator[r]{\sum^M_{k=1}} p_{\left(k\right)}$. Clearly, $\frac{\partial^2 P_{\left(\rm L,avg\right)}}{\partial \Gamma_{\left(1\right)}^2} < 0 $, which implies the objective function of problem $\left(P6\right)$ is concave. Moreover, the subjected constraints $C14$ and $C15$ are linear. Hence, problem $\left(P6\right)$ is a convex optimization problem.
\end{proof}

We can apply the KKT conditions to obtain the global solution for problem $\left(P6\right)$ because it is a convex problem. Hence, the Lagrangian of problem $\left(P6\right)$ is:
\begin{align}
    \mathcal{L}_2 = &- \left[ \eta_{\rm h} P_{\rm a} - \eta_{\rm h} P_{\rm a} \smashoperator[r]{\sum^M_{k=1}} p_{\left(k\right)} \Bigl( \Gamma_{\left(1\right)} - 2 m_{\rm th} \left(k-1\right) \Bigr)^2 \right] \nonumber \\
    &+ \lambda_4 \Bigl( \Gamma_{\rm lb} - \Gamma_{\left(1\right)} + 2 m_{\rm th} \left(M-1\right) \Bigr) \nonumber \\
    &+ \lambda_5 \Bigl( \Gamma_{\left(1\right)} - \Gamma_{\rm ub} \Bigr),
\end{align}
where $\lambda_4$ and $\lambda_5$ are the Lagrange multipliers associated with $C14$ and $C15$, respectively. The KKT point $\left( \Gamma^\ast_{\left(1\right)}, \lambda^\ast_4, \lambda^\ast_5 \right)$ that gives $P^\ast_{\left(\rm L,avg\right)}$ is obtained by solving the following equations.
\begin{align}
    \begin{split}& \frac{\partial \mathcal{L}_2}{\partial \Gamma_{\left(1\right)}} = 2 \eta_{\rm h} P_{\rm a} \smashoperator[r]{\sum^M_{k=1}} p_{\left(k\right)} \Bigl[ \Gamma_{\left(1\right)} - 2 m_{\rm th} \left(k-1\right) \Bigr] \\ & \hspace{50mm} - \lambda_4 + \lambda_5 = 0,\end{split} \label{eq22} \\
    &\lambda_4 \Bigl( \Gamma_{\rm lb} - \Gamma_{\left(1\right)} + 2 m_{\rm th} \left(M-1\right) \Bigr) = 0, \label{eq23} \\
    &\lambda_5 \Bigl( \Gamma_{\left(1\right)} - \Gamma_{\rm ub} \Bigr) = 0. \label{eq24}
\end{align}

along with $\lambda_4, \lambda_5 \geq 0$. Likewise, \eqref{eq22} is the sub-gradient condition, and \eqref{eq23} and \eqref{eq24} are the complementary slackness conditions. The method we use to determine the KKT point is stated in the following Lemma. 

\begin{lemma} \label{L7}
    The $\Gamma^\ast_{\left(1\right)}$ is obtained from one of the 3 cases, which are case (a): $\lambda_4=\lambda_5=0$, case (b): $\lambda_4= 0, \lambda_5 \neq 0$, and case (c): $\lambda_4 \neq 0, \lambda_5=0$.
\end{lemma}

\begin{proof}
    Generally, the KKT point is determined by setting either or both $\lambda_4$ and $\lambda_5$ as either zero or non-zero values. In our problem, we noted that $\lambda_4 \neq 0$ and $\lambda_5 \neq 0$ only occur when $\Gamma_{\rm ub} = \Gamma_{\rm lb} + 2 m_{\rm th} \left(M-1\right)$. When this condition holds, the solution obtained from this case is the same as we obtained from either case (b) or case (c). Therefore, there are only 3 cases to be considered. 
\end{proof}

We denote $\Gamma^{\left(a\right)}_{\left(1\right)}, \Gamma^{\left(b\right)}_{\left(1\right)}$, and $\Gamma^{\left(c\right)}_{\left(1\right)}$ as the solutions obtained from the 3 cases stated in Lemma \ref{L7}. Specifically, in case (a), $\Gamma^{\left(a\right)}_{\left(1\right)}$ is obtained with \eqref{eq22}. Likewise, $\Gamma^{\left(b\right)}_{\left(1\right)}$ in case (b) is obtained with \eqref{eq23}, whereas $\Gamma^{\left(c\right)}_{\left(1\right)}$ in case (c) is obtained with \eqref{eq24}. Subsequently, the closed-form expression for the solutions are:
\begin{subequations}
\begin{align}
    \Gamma^{\left(a\right)}_{\left(1\right)} &= 2 m_{\rm th} \smashoperator[r]{\sum_{k=1}^M} p_{\left(k\right)} \left(k-1\right), \\
    \Gamma^{\left(b\right)}_{\left(1\right)} &= \Gamma_{\rm ub}, \\
    \Gamma^{\left(c\right)}_{\left(1\right)} &= \Gamma_{\rm lb} + 2 m_{\rm th} \left(M-1\right).
\end{align}    
\end{subequations}

Similarly, since $\left(P6\right)$ is a convex problem, only one KKT point with $\lambda^\ast_4, \lambda^\ast_5 \geq 0$ will satisfy all the constraints, which is selected from $\{ \Gamma^{\left(a\right)}_{\left(1\right)}, \Gamma^{\left(b\right)}_{\left(1\right)}, \Gamma^{\left(c\right)}_{\left(1\right)} \}$. Specifically, we notice that $\lambda_5$ from case (b) is positive when $\Gamma^{\left(b\right)}_{\left(1\right)} < \Gamma^{\left(a\right)}_{\left(1\right)}$. In contrast, $\lambda_4$ from the case (c) is positive when $\Gamma^{\left(c\right)}_{\left(1\right)} > \Gamma^{\left(a\right)}_{\left(1\right)}$. Hence, $\Gamma^\ast_{\left(1\right)}$ is determined by:
\begin{equation}
    \Gamma^\ast_{\left(1\right)} = \max \biggl\{ \Gamma^{\left(c\right)}_{\left(1\right)}, \min \Bigl\{ \Gamma^{\left(a\right)}_{\left(1\right)}, \Gamma^{\left(b\right)}_{\left(1\right)} \Bigr\} \biggr\},
    \label{eq36}
\end{equation}

Knowing the sequence of $\Gamma^\ast_{\left(k\right)}$ is a descending order with $k$ on a number line, and therefore: 
\begin{equation}
    \Gamma^\ast_{\left(k+1\right)} = \Gamma^\ast_{\left(k\right)} - 2 m_{\rm th}, \quad \forall k \in \mathbb{M} \setminus M.
\end{equation}
Subsequently, $P^\ast_{\left(\rm L,avg\right)}$ is determined by: 
\begin{equation}
    P^\ast_{\left(\rm L,avg\right)} = \eta_{\rm h} P_{\rm a} \smashoperator[r]{\sum^M_{i=1}} p_{\left(i\right)} \biggl( 1- \left( \Gamma^\ast_{\left(i\right)} \right)^2 \biggr).
\end{equation}

Till this step, we solved problem $\left(P6\right)$ and were able to obtain the closed-form expression for the global solution with the given of $p_{\left(k\right)}$, $\forall k$. However, we do not know the sequence of all $\Gamma^\ast_{{\rm a}i}$, which means $\mathbf{p}_2$ is unknown. One possible solution to determine $\mathbf{p}_2$ is to consider all permutation sequences of $\mathbf{p}$. Using this approach, we realize that there is a maximum of $M!$ possible sequences, each of which can yield a candidate for $\mathbf{\Gamma^\ast_{\rm a}}$.

Considering all the permutation sequences of the symbols' probability for finding the candidate of $\mathbf{\Gamma^\ast_{\rm a}}$ would require considerable computing time, especially when $M$ is large. Nevertheless, from the solving process in Section~\ref{SecB}, we observe a regular pattern of $\hat{\Gamma}^\ast_{{\rm a}i}$, which inspires an approach to simplify the calculation due to the $M!$ possible sequences of $\mathbf{p}$. Specifically, we note that symbols with higher occurrence probabilities have lower mismatch degrees. Consequently, the optimal reflection coefficient corresponding to higher occurrence probability will be set closer to $0$ than those with lower occurrence probabilities. This implies $\Gamma^\ast_{{\rm a}i}$ is dominated by the symbol probability. Hence, the optimal solution of $\left(P2\right)$ will satisfy the relationship $\abs{\Gamma^\ast_{{\rm a}i}} \leq \abs{\Gamma^\ast_{{\rm a}k}}$ for $p_i \geq p_k$, $ \forall i,k \in \mathbb{M}$ and $i \neq k$.

In this context, once we locate $\Gamma_{\rm a1}$, the remaining $\Gamma_{{\rm a}i}$ will be placed alternatively on both sides of $\Gamma_{\rm a1}$. However, it might happen that either the left or right side is truncated due to the lower or upper bound constraint, and the remaining $\Gamma_{{\rm a}i}$ will then be placed sequentially on the other side. This approach significantly reduces the number of candidate solutions from $M!$ to $2 \left(M-1\right)$. An example of this idea for 4-ASK is provided in the following matrix:
\begin{align*}
\begin{bmatrix}
    \Gamma_{\rm a1} & \Gamma_{\rm a2} & \Gamma_{{\rm a}3} & \Gamma_{{\rm a}4} \\
    \Gamma_{{\rm a}3} & \Gamma_{\rm a1} & \Gamma_{\rm a2} & \Gamma_{{\rm a}4} \\
    \Gamma_{\rm a2} & \Gamma_{\rm a1} & \Gamma_{{\rm a}3} & \Gamma_{{\rm a}4} \\
    \Gamma_{{\rm a}4} & \Gamma_{{\rm a}3} & \Gamma_{\rm a1} & \Gamma_{\rm a2} \\
    \Gamma_{{\rm a}4} & \Gamma_{\rm a2} & \Gamma_{\rm a1} & \Gamma_{{\rm a}3} \\
    \Gamma_{{\rm a}4} & \Gamma_{{\rm a}3} & \Gamma_{\rm a2} & \Gamma_{\rm a1} 
\end{bmatrix}
\Longleftrightarrow
\begin{bmatrix}
    p_1 & p_2 & p_3 & p_4 \\
    p_3 & p_1 & p_2 & p_4 \\
    p_2 & p_1 & p_3 & p_4 \\
    p_4 & p_3 & p_1 & p_2 \\
    p_4 & p_2 & p_1 & p_3 \\
    p_4 & p_3 & p_2 & p_1 
\end{bmatrix}
\end{align*}
where the 6 rows in the left matrix are all the possible sequences of $\left[ \Gamma_{\rm a1}, \Gamma_{\rm a2}, \Gamma_{{\rm a}3}, \Gamma_{{\rm a}4} \right]$ with the corresponding matrix of probabilities $\left[ p_1, p_2, p_3, p_4 \right]$ on the right. Therefore, we will only use the probability sequences from the above matrix to obtain $P^\ast_{\rm L,avg}$ for the 4-ASK modulation scheme. A similar concept can be extended to $M$-ASK, where the method to obtain the probability sequences matrix is given in Algorithm \ref{Algo2}.
 
\begin{algorithm}[!t]
    \caption{Generation of Sequence Matrix $\mathbf{A}$}
    \label{Algo2}
    \begin{algorithmic}[1]
        \Require $M, \mathbf{p}$
        \Ensure $\mathbf{A}$
        \State $N=2 (M-1)$
        \State Set $\mathbf{A}^{N \times M}=
            \begin{bmatrix}
                a_{11} & a_{12} & \cdots & a_{1M} \\
                a_{21} & a_{22} & \cdots & a_{2M} \\
                \vdots & \vdots & \vdots & \vdots \\
                a_{N1} &a_{N2} & \cdots & a_{NM}
            \end{bmatrix}$
        \For{$\omega_1 = 1$ to $M$}
            \State $a=\omega_1$
            \For{$\omega_2=1$ to $2$}
                \State $\beta=0$
                \If{$\omega_1 \neq 1 $ or $\omega_2 \neq 1$}
                    \If{$\omega_1 \neq M$ or $\omega_2 \neq 2$}
                        \For{$\omega_3=1$ to $M$}
                            \If{$\omega_2=1$}
                                \State $a \leftarrow a+(-1)^{\omega_3+1} (\omega_3-1)$
                            \ElsIf{$\omega_2=2$}
                                \State $a \leftarrow a+(-1)^{\omega_3} (\omega_3-1)$
                            \EndIf
                            \State $n=2\omega_1-( \frac{1}{2}+\frac{1}{2}(-1)^{\omega_2+1} )-1$
                            \If{$a \leq M$ and $a \geq1$ and $\beta=0$}
                                \State $m=a$
                            \Else 
                                \If{$\omega_1 \leq \frac{M}{2}$}
                                    \State $\beta=\omega_3$
                                \Else
                                    \State $\beta=M-\omega_3+1$
                                \EndIf
                                \State $m=\beta$
                            \EndIf
                            \State $a_{nm} = p_{\omega_3}$
                        \EndFor
                    \EndIf
                \EndIf
            \EndFor
        \EndFor
        \State Return $\mathbf{A}$
    \end{algorithmic}
\end{algorithm}

\begin{algorithm}[ht]
    \caption{Optimal Reflection Coefficient for $M$-ASK Modulation}
    \label{Algo1}
    \begin{algorithmic}[1]
    \Require $M$, $m_{\rm th}$, $P_{\rm sen}$, $P_{\rm L,min}$, $P_{\rm b,min}$, $P_{\rm t}$, $G_{\rm t}$, $G_{\rm r}$, $d$, $f$, $n$, $ \eta_{\rm h}$ and $ \eta_{\rm b}$
    \Ensure $\mathbf{\Gamma_a}, P^\ast_{\rm L,avg}$
    \State Run Algorithm \ref{Algo2} to obtain $\mathbf{A}$
    \State Set $\mathbb{T}=0$
    \For{$\nu=1$ to $2(M-1)$}
        \For{$u=1$ to $M$}
            \State $p_{\left(u\right)} = a_{\nu u}$
        \EndFor
        \State $\Gamma^\ast_{\left(1\right)} = \max \biggl\{ \Gamma^{\left(c\right)}_{\left(1\right)}, \min \Bigl\{ \Gamma^{\left(a\right)}_{\left(1\right)}, \Gamma^{\left(b\right)}_{\left(1\right)} \Bigr\} \biggr\}$
        \For{$w=1$ to $M$}
            \State $\Gamma^\ast_{\left(w\right)} = \Gamma^\ast_{\left(1\right)} - 2 m_{\rm th} \left(w-1\right)$
        \EndFor
        \State $P^\ast_{\left(\rm L,avg\right)} = \eta_{\rm h} P_{\rm a} \smashoperator[r]{\sum^M_{k=1}} p_{\left(k\right)} \biggl( 1- \left( \Gamma^\ast_{\left(k\right)} \right)^2 \biggr)$
        \If{$P^\ast_{\left(\rm L,avg\right)} > \mathbb{T}$}
            \State $\mathbb{T} = P^\ast_{\left(\rm L,avg\right)}$
            \State $\mathbf{B}^{2 \times M} = $
                $\begin{bmatrix}
                    p_{\left(1\right)} & p_{\left(2\right)} & ... & p_{\left(M\right)}\\
                    \Gamma^\ast_{\left(1\right)} & \Gamma^\ast_{\left(2\right)} & ... & \Gamma^\ast_{\left(M\right)} 
                \end{bmatrix}$
        \EndIf
    \EndFor
    \State $\mathbf{C}^{2 \times M} = sort \left( \mathbf{B}, 1 \right)$
    \For{$i=1$ to $M$}
        \State $\Gamma^\ast_{{\rm a}i} = c_{2i}$
    \EndFor
    \State Return $\mathbf{\Gamma_a} = \left[ \Gamma^\ast_{{\rm a}1}, \Gamma^\ast_{{\rm a}2},...\Gamma^\ast_{{\rm a}M} \right]$, $P^\ast_{\rm L,avg} = \mathbb{T}$
    \end{algorithmic}
\end{algorithm}

Algorithm 1 provides us with matrix $\mathbf{A}$, and from there, we use each row of $\mathbf{A}$ as $\mathbf{p}_2$ to obtain all the candidate solutions with~\eqref{eq36}. During this process, we make use of the variable $\mathbb{T}$ and matrix $\mathbf{B}$, to keep track of $P^\ast_{\left(\rm L,avg\right)}$ and the corresponding $p_{\left(k\right)}$ and $\Gamma^\ast_{\left(k\right)}$, $\forall k \in \mathbb{M}$. As we iterate through each probability sequence in $\mathbf{A}$, we update $\mathbb{T}$ and $\mathbf{B}$ whenever $P^\ast_{\left(\rm L,avg\right)}$ is greater than the current $\mathbb{T}$. This process continues until all probability sequences in $\mathbf{A}$ have been used, and we will obtain the optimal $\mathbf{p}_2$ and $\mathbf{\Gamma}_{(a)}$ that yields $P^\ast_{\rm L,avg}$. Given that $\mathbf{p}_2$ represents one of the rows in matrix $\mathbf{A}$, and $p_i$ in $\mathbf{p}$ is arranged in descending order, we can obtain $\Gamma^\ast_{{\rm a}i}$ from $\mathbf{\Gamma_{(a)}}$ by rearranging $\Gamma_{(k)}$ in $\mathbf{\Gamma_{(a)}}$ such that the corresponding $p_{(k)}$ in $\mathbf{p}_2$ are arranged in descending order. By doing so, we can determine the global solution $\mathbf{\Gamma^\ast_{\rm a}}$ of problem $\left(P2\right)$ and $P^\ast_{\rm L,avg}$. These steps are summarized in Algorithm \ref{Algo1}. The $sort(\mathbf{E},x)$ function, in Algorithm \ref{Algo1}, sorts the columns of the input matrix $\mathbf{E}$ in descending order based on the elements of the $x$ row.\footnote{It should be noted that our proposed method for maximizing the average harvested power of passive backscatter tags can be easily adapted to solve any nonlinear energy harvesting model, provided it is a monotonic increasing function of input power. An example of applying this method to a nonlinear model is detailed in Appendix~\ref{ApdC}.}

\section{Results and Discussion} \label{sectionRD}

 We numerically demonstrate the tag performance with our proposed optimal reflection coefficients. Specifically, we determine and observe the achievable $P^\ast_{\rm L,avg}$ for different tag design parameters and BackCom systems. Unless otherwise stated, we set $P_{\rm t}=1$ W, $f=915$ MHz, $c=3 \times 10^8$ ms$^{-1}$, $G_{\rm t}=4$, $G_{\rm r}=1.5$, $d_{\rm o} = 1$ m, $d=8$ m, $P_{\rm sen} = 2$ $\mu$W, $P_{\rm L,min}=5$ $\mu$W, $P_{\rm b,min}=3$ $\mu$W, $\eta_{\rm h} = 0.4$, $\eta_{\rm b} = 0.8$, $\sigma^2 = -90$ dBm and $R_{\rm r}=50~\Omega$. Besides, we consider our system model operates in obstructed factories, where the path loss exponent is set to $n=2.5$~\cite{miranda2013path}. Depending on the modulation order $M$, each symbol's total number of bits will be $\log_2 M$. A sequence of bits ‘1’ and ‘0’ with the probability $\mathcal{P}^{\{1\}}$ and $\mathcal{P}^{\{0\}}$, respectively are used to generate the symbols. For example, considering 4-ASK, the transmit symbols are ‘11’, ‘10’, ‘01’, and ‘00’. The corresponding symbol probabilities are $p_1=\left( \mathcal{P}^{\{1\}} \right)^2$, $p_2= \mathcal{P}^{\{1\}} \mathcal{P}^{\{0\}}$, $p_3= \mathcal{P}^{\{0\}} \mathcal{P}^{\{1\}}$, and $p_4=\left( \mathcal{P}^{\{0\}} \right)^2$. Without the loss of generality, we set $\mathcal{P}^{\{1\}} \geq \mathcal{P}^{\{0\}}$ to simplify our numerical investigation, which can be easily switched to fit the opposite case. 
 
\subsection{Insight on Optimal Reflection Coefficients}

\begin{figure}[!t]
\centering
\includegraphics[width=3.48in]{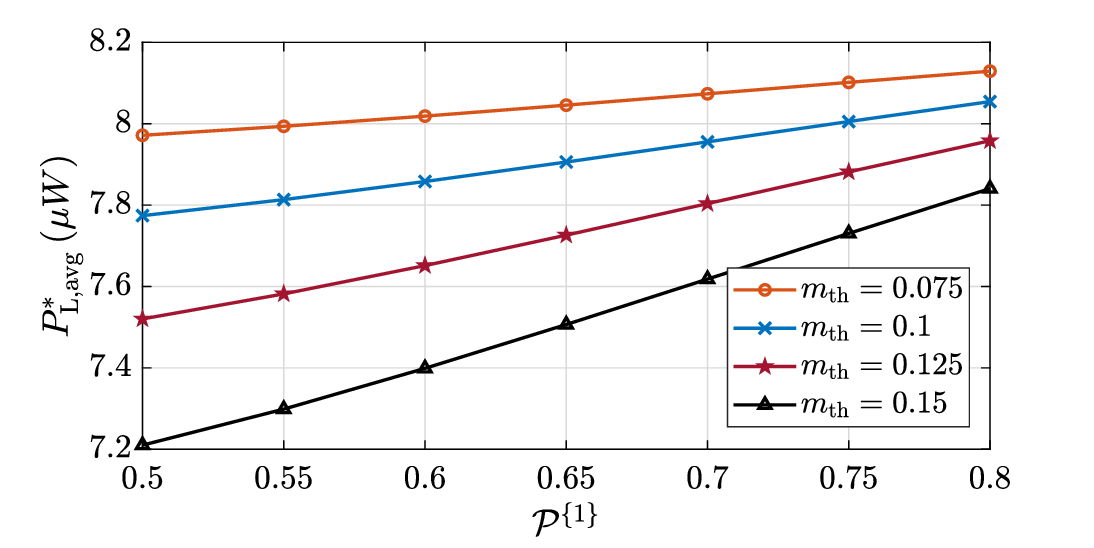}
\caption{Design insights on $P^\ast_{\rm L,avg}$ for 4-ASK.}
\label{figa222}
\end{figure}

\begin{figure}[!t]
\centering
\includegraphics[width=3.48in]{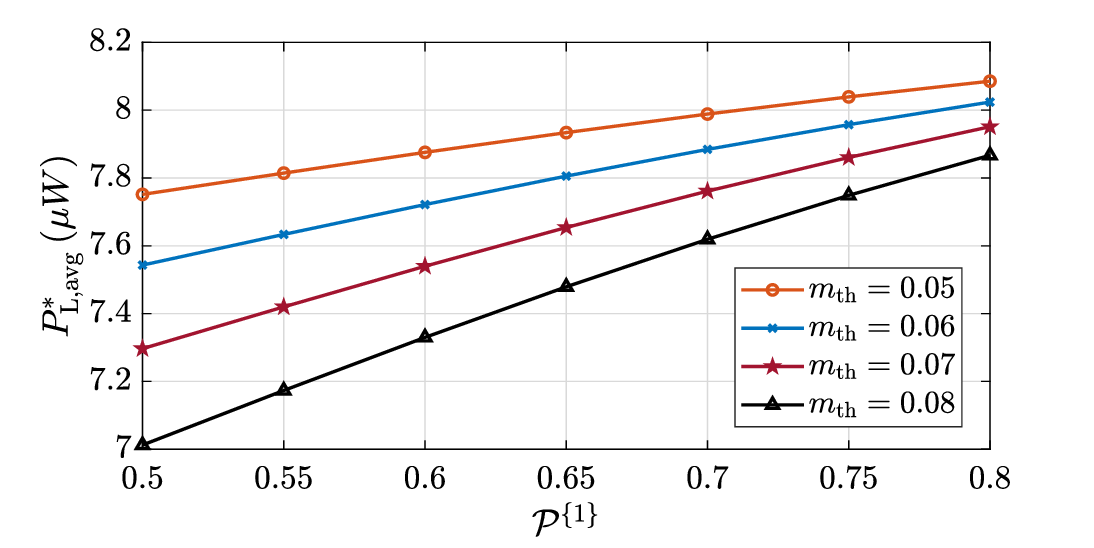}
\caption{Design insights on $P^\ast_{\rm L,avg}$ for 8-ASK.}
\label{figb222}
\end{figure}

\begin{figure}[!t]
    \centering
    \subfigure[]
    {\hspace{-0.4mm}\includegraphics[width=1.8in]{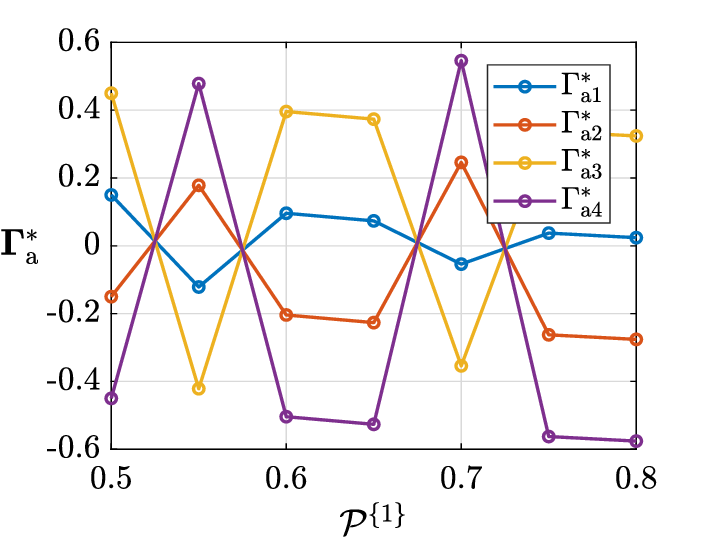}
        \label{fig3aaa}
        \hspace{-0.175in}
    }
    \subfigure[]
    {
        \hspace{-0.175in}
        \includegraphics[width=1.8in]{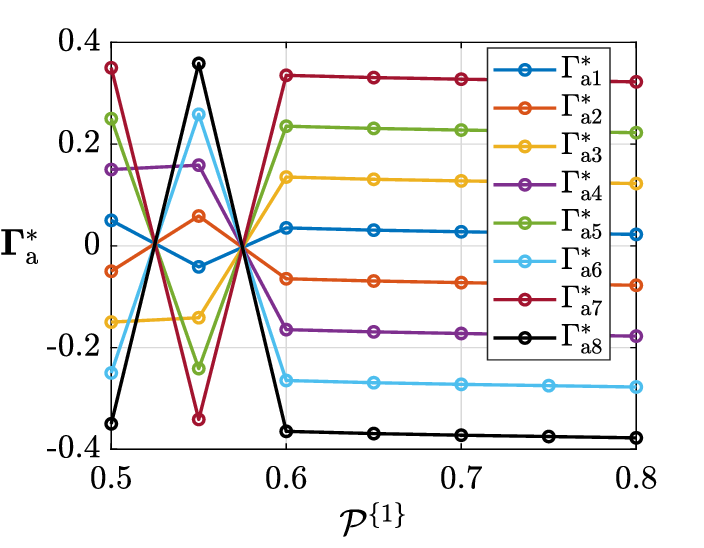}
        \label{fig3bbb}
        \hspace{-0.2in}
    }
    \caption{$\Gamma^\ast_{{\rm a}i}$ versus $\mathcal{P}^{\{1\}}$ at $d=8$m with (a) 4-ASK, $m_{\rm th}=0.15$ and (b) 8-ASK, $m_{\rm th}=0.05$.}
    \label{fig3ccc}
\end{figure}

In practical scenarios, the sensor data and the tag's unique identifier may exhibit a significant bias towards either bit '1' or '0'. We leverage the non-equiprobable occurrence probability of symbols to maximize the tag's harvested power, as it is highly dependent on the reflection coefficient associated with the probability of the transmitted symbols. We determine the correlation between $P^\ast_{\rm L,avg}$ and $\mathcal{P}^{\{1\}}$ at different $m_{\rm th}$ for 4-ASK and 8-ASK modulation schemes, which are shown in Fig.~\ref{figa222} and~\ref{figb222}, respectively. It is crucial to note that the selection of $m_{\rm th}$ is guided by ensuring the tag meets the operational requirements of the BackCom system. We observe that $P^\ast_{\rm L,avg}$ increases with $\mathcal{P}^{\{1\}}$, which implies the significance of considering the transmit symbols' probability in the tag design. To understand how it improves $P^\ast_{\rm L,avg}$, we study the underlying optimal reflection coefficient $\Gamma^\ast_{{\rm a}i}$ of each state. Fig. \ref{fig3aaa} and \ref{fig3bbb} illustrate $\Gamma^\ast_{{\rm a}i}$ versus $\mathcal{P}^{\{1\}}$ with different $m_{\rm th}$ for 4-ASK and 8-ASK, respectively. The $\Gamma^\ast_{{\rm a}i}$ corresponding to the symbol with a relatively higher occurrence probability is closer to 0, resulting in higher $P_{{\rm L}i}$ and thereby increasing $P^\ast_{\rm L,avg}$. 

\begin{figure}[!t]
    \centering
    \includegraphics[width=3.48in]{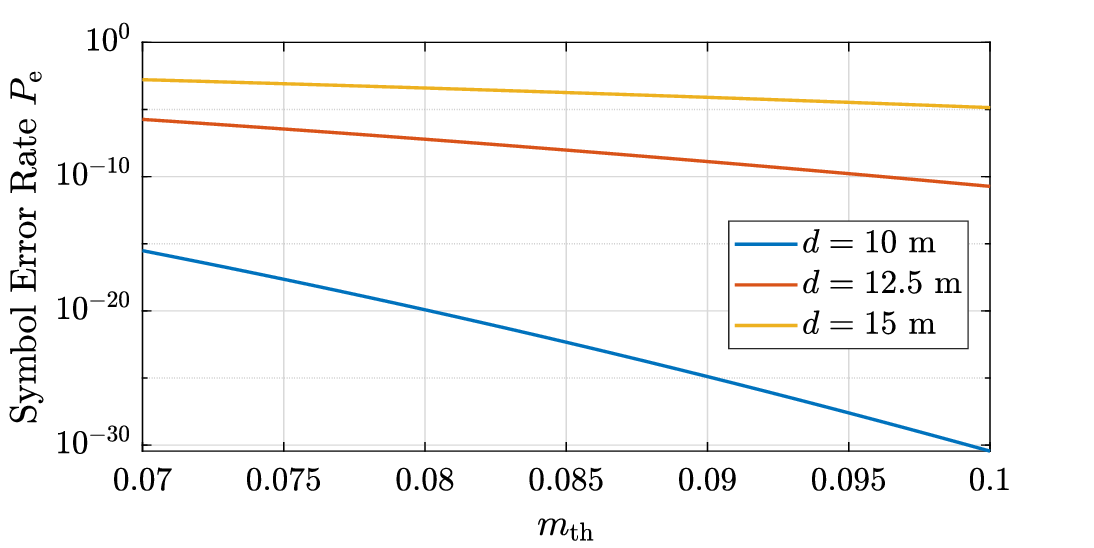}
    \caption{Symbol error rate of BackCom system versus modulation index $m_{\rm th}$ for different transmission distances.}
    \label{figBER}
\end{figure}

On the other hand, the BackCom system's SER performance depends on $m_{\rm th}$ as given in \eqref{eq5}. To illustrate this relationship, we plot the SER for various values of $m_{\rm th}$ and different transmission distances in Fig.~\ref{figBER}, visually demonstrating the impact of $m_{\rm th}$ on the SER. It is evident that the SER decreases with higher $m_{\rm th}$, resulting in reliable BackCom. In Fig. \ref{figa222} and \ref{figb222}, we observed a decrease in $P^\ast_{\rm L,avg}$ with an increase in $m_{\rm th}$. This is because it requires a greater mismatch degree of $\Gamma^\ast_{{\rm a}i}$ to achieve the high $m_{\rm th}$, thereby reducing $P_{{\rm L}i}$ in each state and drops $P^\ast_{\rm L,avg}$ in the sense.

\begin{figure}[!t]
    \centering
    \includegraphics[width=3.6in]{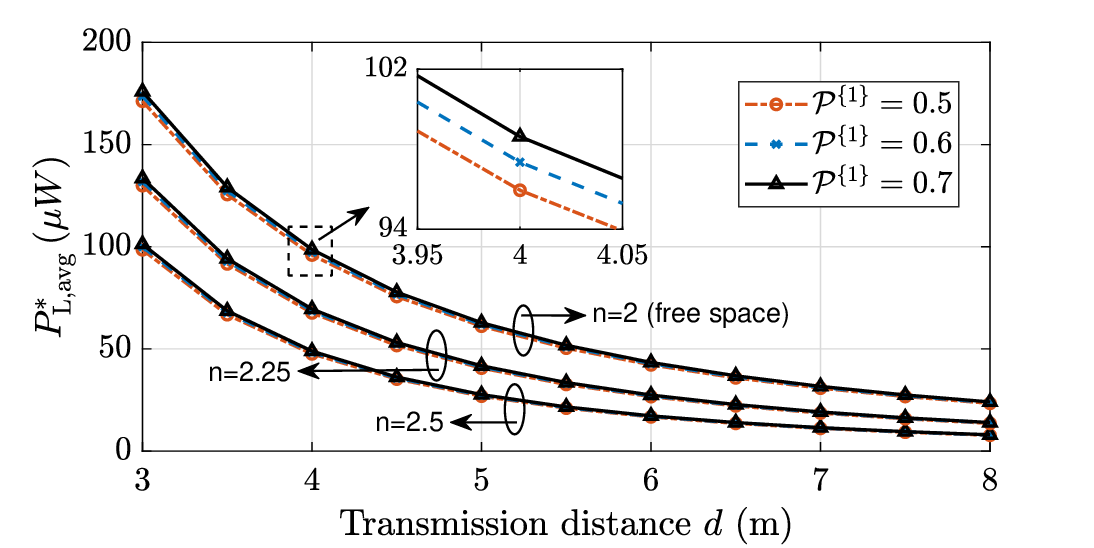}
    \caption{$P^\ast_{\rm L,avg}$ versus $d$ at different $\mathcal{P}^{\{1\}}$ for 8-ASK.}
    \label{fig1122}
\end{figure}

\begin{figure}[!t]
    \centering
    \subfigure[]
    {\hspace{-0.4mm}\includegraphics[width=1.8in]{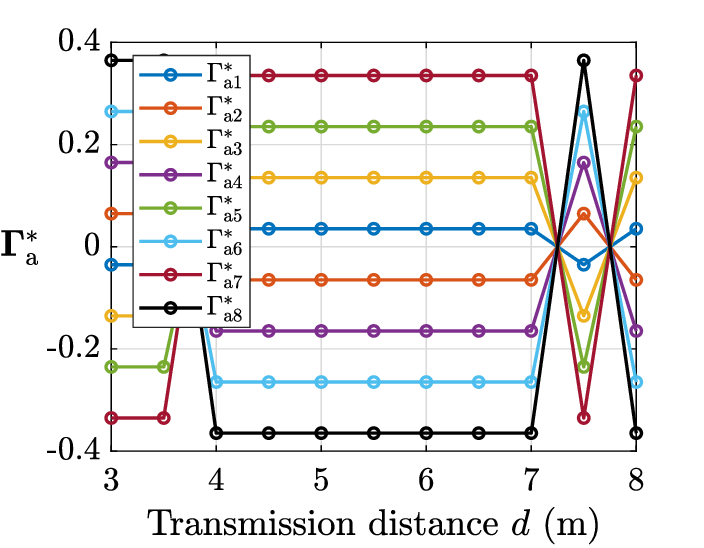}
        \label{fig321a}
        \hspace{-0.175in}
    }
    \subfigure[]
    {
        \hspace{-0.175in}
        \includegraphics[width=1.8in]{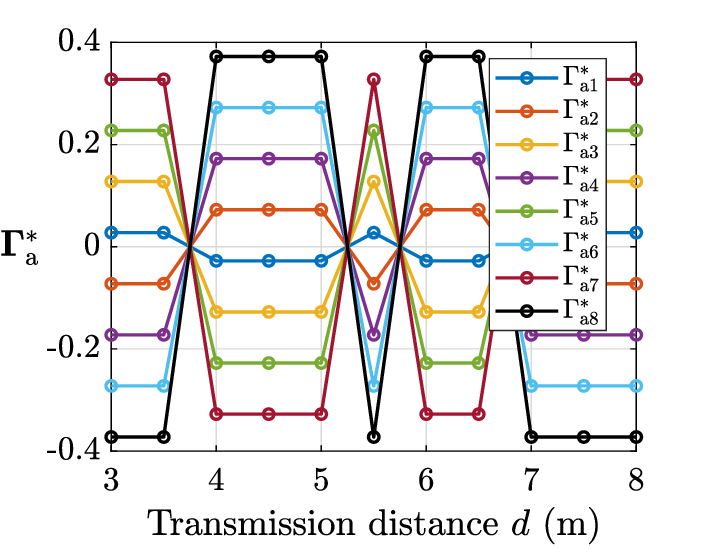}
        \label{fig321b}
        \hspace{-0.2in}
    }
    \caption{Optimal Reflection coefficient for 8-ASK at (a) $\mathcal{P}^{\{1\}}=0.6$ and (b) $\mathcal{P}^{\{1\}}=0.7$ for $n=3$.}
    \label{fig321}
\end{figure}

As the passive tag relies on harvesting energy from the reader's broadcasted signal, the transmission distance between the reader and the tag emerges as a crucial factor influencing the tag's harvesting power. Generally, the transmitted signal's energy experiences a significant decrease over propagation distance. In consideration of this, we set $m_{\rm th}=0.05$, $P_{\rm b,min} = 5 \mu$W, and $M=8$ to plot $P^\ast_{\rm L,avg}$ against $d$, as illustrated in Fig.~\ref{fig1122}. Additionally, we show the impact of the path loss exponent $n$ on $P^\ast_{\rm L,avg}$ within the same plot.

Results in Fig.~\ref{fig1122} demonstrate a substantial decrease in $P^\ast_{\rm L,avg}$ as $d$ increases, dropping to less than $6.3\mu$W when $d=8$m and $n=3$. We also observe a significant decrease in $P^\ast_{\rm L,avg}$ as $n$ increases from 2 to 2.5 and subsequently to 3. This underscores that the BackCom system environment is vital in energy harvesting. The decrease in $P^\ast_{\rm L,avg}$ can be attributed to the exponential attenuation of the RF carrier over the distance $d$ and the presence of a poor propagation medium, leading to a significant reduction in $P_{\rm a}$. Besides, the $d$ and $n$ exhibit a more pronounced impact on $P^\ast_{\rm L,avg}$ compared to $\mathcal{P}^{\{1\}}$. Similarly, we plot the underlying $\mathbf{\Gamma^\ast_{\rm a}}$ to observe how each state's reflection coefficient varies with $d$. Specifically, Fig.~\ref{fig321a} and~\ref{fig321b} depict plots for $\mathcal{P}^{\{1\}}=0.6$ and $\mathcal{P}^{\{1\}}=0.7$, respectively. Our observations reveal that $\abs{\Gamma^\ast_{{\rm a}i}}$ remains constant within a certain transmission range (up to $d=7.5$m in our case) and only varies beyond a specific distance ($d=8$m in our case). This suggests that $d$ does not influence $\Gamma^\ast_{{\rm a}i}$ within a certain range where none of the $\Gamma^\ast_{{\rm a}i}$ values are selected that just satisfy $\Delta P_{\rm L,min}$ or $P_{\rm b,min}$. As $d$ increases and results in a decrease in $P_{\rm a}$, either or both $\Delta P_{\rm L,min}$ and $P_{\rm b,min}$ constraints set a boundary to $\Gamma^\ast_{{\rm a}i}$. Specifically, in Fig.~\ref{fig321b}, we observe that the $P_{\rm b,min}$ constraint establishes a boundary such that $\Gamma^\ast_{{\rm a}i} \leq 0.277$ when $d=8$m. This underscores the significance of our optimal reflection coefficient selection algorithm in maximizing the average harvested power while ensuring the essential operational requirements for the ASK-modulated BackCom system are met.

\subsection{Performance Comparison}
\subsubsection{General Case for $M$-ASK}

\begin{figure}
    \centering
    \includegraphics[width=3.3in]{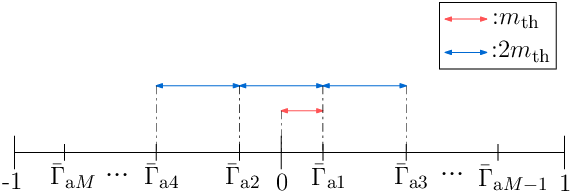}
    \caption{Distribution of $\bar{\Gamma}_{{\rm a}i}$ on a number line.}
    \label{figNum}
\end{figure}

In this study, we aim to compare our load selection method with a conventional approach, focusing on demonstrating the enhancement in $P_{\rm L,avg}$ Recent years have seen various reflection coefficient selections for backscatter tag design, with studies such as~\cite{nikitin2005power,nikitin2007differential,bletsas2010improving} considering an equal mismatch of the reflection coefficient for the BASK modulation scheme. Adopting a similar concept, we introduce a reflection coefficient selection, denoted as $\bar{\Gamma}_{{\rm a}i}$, $\forall i \in \mathbb{M}$. Specifically, $\bar{\Gamma}_{{\rm a}i}$ is symmetrically placed at $0$ in the descending order of $p_i$, with an equal separation of $2 m_{\rm th}$. Assuming $\bar{\Gamma}_{a1} > \bar{\Gamma}_{a2}$ without loss of generality, we define $\bar{\Gamma}_{{\rm a}i} = m_{\rm th} \Biggl[ 1 - \left(-1\right)^i \biggl( i + \frac{1}{2} \Bigl( \left(-1\right)^i - 1 \Bigr) \biggl) \Biggr]$. The symmetrical distribution of $\bar{\Gamma}_{{\rm a}i}$ on a number line is illustrated in Fig.~\ref{figNum}, visualizing the mismatch degree of $\bar{\Gamma}_{{\rm a}i}$ at each state. Accordingly, we define $\bar{P}_{\rm L,avg}$ as the average harvested power determined by $\bar{\Gamma}_{{\rm a}i}$ as the benchmark to highlight the merits of our optimal reflection coefficients selection. 

\begin{figure}[!t]
\centering
\includegraphics[width=3.48in]{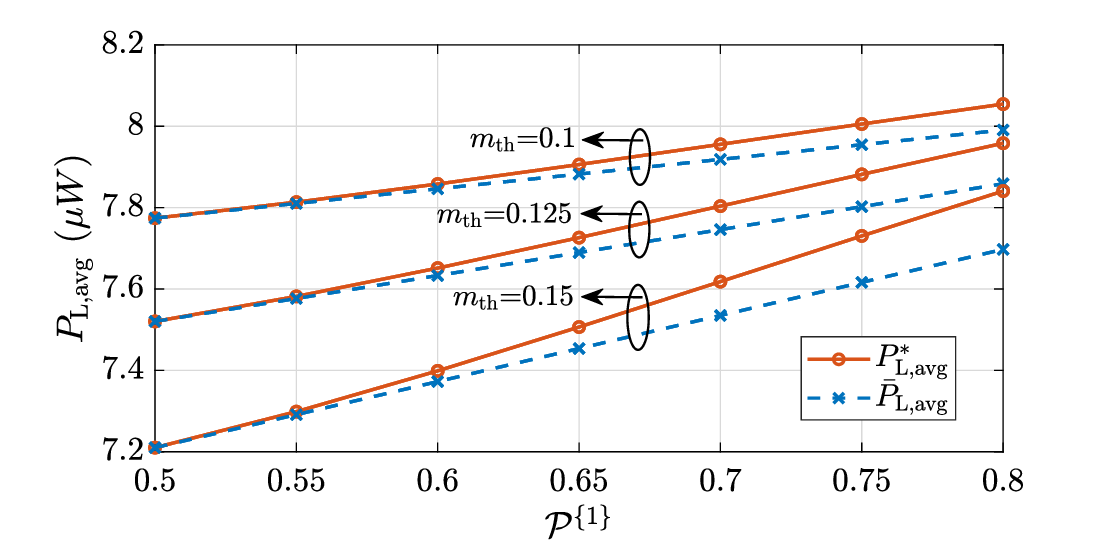}
\caption{$P_{\rm L,avg}$ at different $\mathcal{P}^{\{1\}}$ for 4-ASK.}
\label{figCmp}
\end{figure}

Fig.~\ref{figCmp} plots $P^\ast_{\rm L,avg}$ and $\bar{P}_{\rm L,avg}$ for varying $\mathcal{P}^{\{1\}}$ and different $m_{\rm th}$ at $d=8$m and $M=4$. The result shows that $P^\ast_{\rm L,avg}$ is always greater than $\bar{P}_{\rm L,avg}$, with the relative improvement increasing with $\mathcal{P}^{\{1\}}$. This suggests that our optimal solution for tag design exhibits superior performance. Moreover, the symmetrical arrangement of $\bar{\Gamma}_{{\rm a}i}$ imposes limitations on its application because of the asymmetrical boundary constraints of $\Gamma_{{\rm a}i}$. Conversely, our algorithm offers a more adaptable and flexible reflection coefficient selection, particularly allowing a higher $m_{\rm th}$.

\subsubsection{The Special Case with BASK}\label{B}

\begin{figure}[!t]
	\centering
    \includegraphics[width=3.48in]{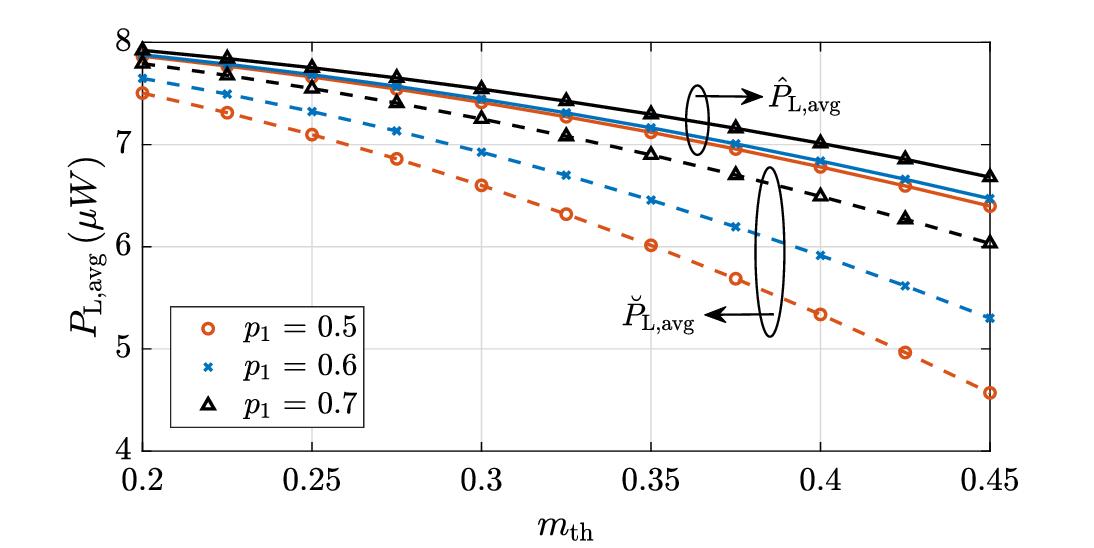}
	\caption{$P_{\rm L,avg}$ in BASK at different $m_{\rm th}$.}
	\label{Fig11abc}
\end{figure}

\begin{figure}[!t]
	\centering
    \includegraphics[width=3.48in]{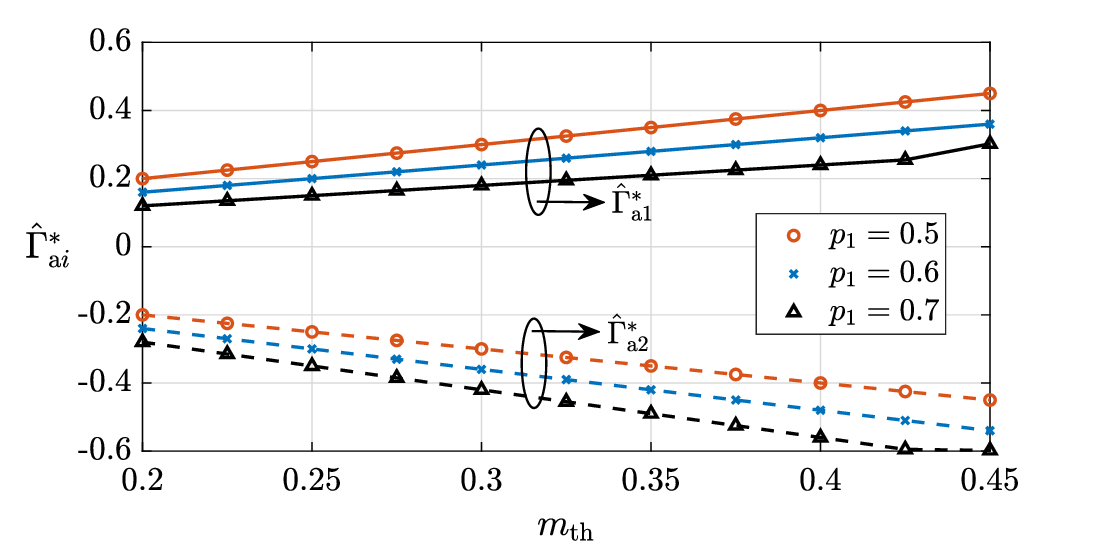}
	\caption{$\hat{\Gamma}^\ast_{{\rm a},i}$ versus $m_{\rm th}$.}
	\label{Fig12abc}
\end{figure}

 As mentioned, BASK is the most popular and widely used modulation in the current backscatter tag industry, particularly for RFID tags. In this context, we study the performance of our load selection with the BASK modulation scheme specified in \ref{IV.B}. Additionally, we consider another prevalent load selection in BASK modulation. This involves one load in perfectly matched condition and another greatly mismatched to meet the $m_{\rm th}$ requirement. We denote it as $\breve{\Gamma}_{{\rm a}i}$ and set $\breve{\Gamma}_{{\rm a}1} > \breve{\Gamma}_{{\rm a}2}$. Hence, $\breve{\Gamma}_{{\rm a}1} = 0$ and $\breve{\Gamma}_{{\rm a}2} = - 2 m_{\rm th}$, which are used to determined the average harvested power $\breve{P}_{\rm L,avg}$. 
 
 Subsequently, we set $\breve{P}_{\rm L,avg}$ as the benchmark and graphically compare it with $\hat{P}^\ast_{\rm L,avg}$, as depicted in Fig.~\ref{Fig11abc}. Similarly, we observed that $P_{\rm L,avg}$ decreases with $m_{\rm th}$. More importantly, it is clear that $\hat{P}^\ast_{\rm L,avg}$ is always greater than $\breve{P}_{\rm L,avg}$, which signifies the utility of our design. The average gain achieved by the proposed optimal reflection coefficients over the benchmark is $9.65\%$. Fig.~\ref{Fig12abc} provides insights into the underlying $\hat{\Gamma}^\ast_{{\rm a}i}$, showing that $\hat{\Gamma}^\ast_{{\rm a}2}$ reaches a constant as $m_{\rm th}$ increases to ensure $\hat{\Gamma}^\ast_{{\rm a}2}$ satisfy the $\Delta P_{\rm L,min}$ constraint.

\begin{figure}[t]
    \centering
    \includegraphics[width=\linewidth]{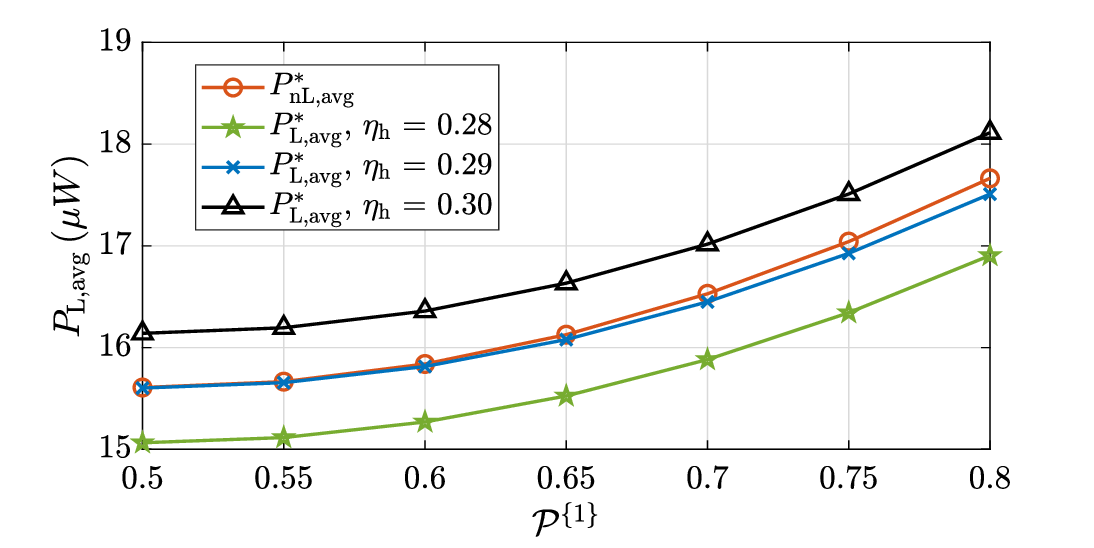}
    \caption{Comparison of maximum harvested power for linear and nonlinear models at different $\mathcal{P}^{\{1\}}$ for the BASK modulation scheme at $d = 5$ m and $m_{\rm th} = 0.5$.}
    \label{figPLNL}
\end{figure}

\begin{figure}[t]
    \centering
    \includegraphics[width=\linewidth]{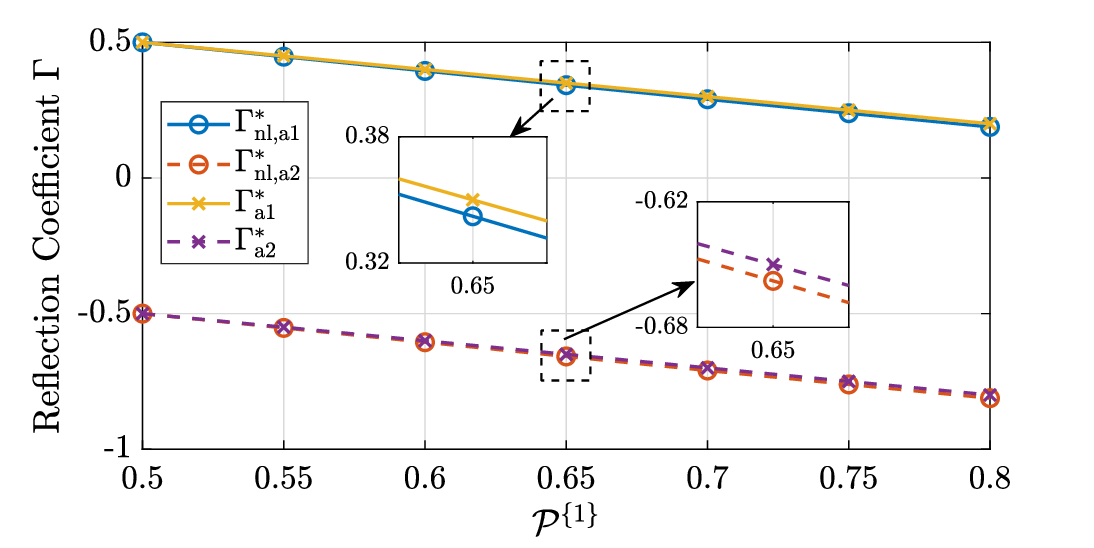}
    \caption{The underlying reflection coefficients $\Gamma^\ast_{{\rm nl,a}i}$ and $\Gamma^\ast_{{\rm a}i}$ that yield $P^\ast_{\rm nL,avg}$ and $P^\ast_{\rm L,avg}$ in Fig. 14, respectively. The $\Gamma^\ast_{\rm a1}$ remains the same for all value of $\eta_{\rm h} = 0.28, 0.29, 0.3$.}
    \label{figRCLNL}
\end{figure}

\subsection{Comparison of CLC and Nonlinear Energy Harvesting Models}

Although we assume the harvested power at the tag follows a CLC model, this process generally exhibits nonlinear behavior due to the implementation of one or more diodes~\cite{alevizos2018sensitive,mishra2017utility,clerckx2019fundamentals}. Therefore, we are interested in comparing the performance of our considered CLC model with the saturation nonlinear model as studied in~\cite{clerckx2019fundamentals}. The maximization of average harvested power (denoted as $P_{\rm nL,avg}$) based on the saturation nonlinear model is detailed in Appendix~\ref{ApdC}. For this comparison, we set $P_{\rm sat} = 0.1071$ mW, $\varkappa = 8000$, and $\varpi = 0.0001$ to determine the optimal reflection coefficient $\Gamma^\ast_{{\rm nl,a}i}$ that results in the maximum harvested power $P^\ast_{\rm nL,avg}$ of the nonlinear model.

Fig.~\ref{figPLNL} demonstrates the maximum harvested power for both CLC and nonlinear models across various probabilities of bit ‘1’, with different energy harvesting efficiencies for the CLC model. Clearly, we observe that $P^\ast_{\rm L,avg}$ increases with greater $\eta_{\rm h}$, as a larger $\eta_{\rm h}$ allows the tag to harvest more power with the same input power. Most importantly, we observe that the difference between $P^\ast_{\rm nL,avg}$ and $P^\ast_{\rm L,avg}$ with $\eta_{\rm h} = 0.29$ is minimal,  with the overall percentage difference being $0.369\%$. Furthermore, the optimal reflection coefficients $\Gamma^\ast_{{\rm nl,a}i}$ and $\Gamma^\ast_{{\rm a}i}$ exhibits a small difference, as shown in Fig.~\ref{figRCLNL}. This implies that the saturation nonlinear energy harvesting model can be accurately approximated by the CLC model, which is proportional to the input power, by considering a constant energy harvesting efficiency with a small discrepancy. It should be noted that approximating the actual energy harvesting model with the CLC model is only valid when the model is a monotonic increasing function of input power, which is the basis of our solution methodology.

\section{Conclusion}\label{sec:conclusion}
 This paper provides a novel load selection for the passive backscatter tag design with the $M$-ASK modulation scheme that maximizes the average harvested power while meeting the tag power sensitivity, reader sensitivity, and SER requirement. Although all the transmit symbols are generally assumed to have equal occurrence probability, we consider them unequal and are application-dependent constants in our work. We studied 2 different tag designs, where the first design is a general $M$-ASK modulator, and the second is a BASK modulator. Besides, we considered more comprehensive operational requirements in the former design. In contrast, the latter design neglected reader sensitivity, which is a common assumption in conventional RFID applications. We obtain global solutions in the closed-form expression with the proposed algorithms. The simulation results have shown that the symbol probability and modulation index can significantly impact the maximum average harvested power. Besides, we found that the average harvested power with the optimal load selection provided a significant average gain over the benchmark, which signifies our proposed reflection coefficient selection. By comparing the CLC model with the saturation nonlinear model, we determined that the CLC model can practically approximate the actual energy harvesting behavior, provided the model is monotonically increasing with respect to input power. This improvement further validates the robustness and applicability of our proposed methodology in real-world scenarios.


\vspace{-2mm}
 
\begin{appendices}
\setcounter{equation}{0}
\setcounter{figure}{0}
\renewcommand{\theequation}{A.\arabic{equation}}
\renewcommand{\thefigure}{A.\arabic{figure}}
 \section{Proof of Lemma 2}\label{ApdA}
Here, we aim to prove Lemma 2 by comparing the backscattered power for different load selections, given the same harvested power. We examine 2 cases, where the first case assumed the reflection coefficient $\Gamma^{(1)} = \Gamma^{(1)}_{\rm a} + j \Gamma^{(1)}_{\rm b}$ with $\Gamma^{(1)}_{\rm b} = 0$.
The second case assumed the reflection coefficient $\Gamma^{(2)} = \Gamma^{(2)}_{\rm a} + j \Gamma^{(2)}_{\rm b}$, with $\Gamma^{(2)}_{\rm a} = 0$. The harvested power both cases is given by:
\begin{align}
    P^{(1)}_{\rm L} &= \eta_{\rm h} \left[ P_{\rm a} \left( 1 - \left( \Gamma^{(1)}_{\rm a} \right)^2 \right) - P_{\rm sen} \right] , \label{eqd111}\\
    P^{(2)}_{\rm L} &= \eta_{\rm h} \left[ P_{\rm a} \left( 1 - \left( \Gamma^{(2)}_{\rm b} \right)^2 \right) - P_{\rm sen} \right], \label{eqd222}
\end{align}
and the backscattered power of these 2 cases are:
\begin{align}
    P^{(1)}_{\rm b} &= \eta_{\rm b} P_{\rm a} G_{\rm t} \left( 1 - \Gamma^{(1)}_{\rm a} \right)^2 \nonumber \\
            &\stackrel{(s_4)}{=} \eta_{\rm b} P_{\rm a} G_{\rm t} \left( 2 - \mathcal{G} + 2 \sqrt{1 - \mathcal{G}} \right) \label{eq17}, \\
    P^{(2)}_{\rm b} &= \eta_{\rm b} P_{\rm a} G_{\rm t} \left( 1 + \left( \Gamma^{(2)}_{\rm b} \right)^2 \right) \nonumber \\
            &\stackrel{(s_5)}{=} \eta_{\rm b} P_{\rm a} G_{\rm t} \left( 2 - \frac{P^{(2)}_{\rm L} + \eta_{\rm h} P_{\rm sen} }{ \eta_{\rm h} P_{\rm a}} \right),
\end{align}
where $(s_4)$ and $(s_5)$ are obtained using \eqref{eqd111} and \eqref{eqd222}, respectively, and $\mathcal{G} = \frac{P^{(1)}_{\rm L} + \eta_{\rm h} P_{\rm sen} }{ \eta_{\rm h} P_{\rm a}}$. By selecting load impedance such that $P^{(1)}_{\rm L} = P^{(2)}_{\rm L}$, it is evident that $P^{(1)}_{\rm b} \geq P^{(2)}_{\rm b} $. Hence, we proved Lemma 2.

\vspace{-2mm}

\section{Proof of Lemma \ref{L4}}\label{ApdB}
We determine the Hessian matrix of the objective function in problem $\left(P4\right)$, which is given as: 
\begin{align}
    \mathbb{H} &=
    \begin{bmatrix}
    \frac{\partial^2 \hat{P}_{\rm L,avg}}{\partial \Gamma_{\rm a1}^2}  & \frac{\partial^2 \hat{P}_{\rm L,avg}}{\partial \Gamma_{\rm a1} \partial \Gamma_{\rm a2}} \\
    \frac{\partial^2 \hat{P}_{\rm L,avg}}{\partial \Gamma_{\rm a2} \partial \Gamma_{\rm a1}} & \frac{\partial^2 \hat{P}_{\rm L,avg}}{\partial \Gamma_{\rm a2}^2}
    \end{bmatrix} \nonumber \\
    &= 
    \begin{bmatrix}
     -2 p_1 \eta_{\rm h} P_{\rm a} & 0 \\ 
    0 & -2 \left(1-p_1\right) \eta_{\rm h} P_{\rm a}    
    \end{bmatrix}.
\end{align}
 Then, we find the eigenvalues of matrix $\mathbb{H}$, which are: 
 \begin{align}
     \lambda_{\rm e1} &= -2 p_1 \eta_{\rm h} P_{\rm a} \hspace{1mm}, \\
     \lambda_{\rm e2} &= -2 \left(1-p_1\right) \eta_{\rm h} P_{\rm a}.
 \end{align}
 
 We observed that the diagonal entries of $\mathbb{H}$ are $\leq 0$, and the determinant of $\mathbb{H}$ being non-negative, $\abs{\mathbb{H}} \geq 0$. Moreover, both the eigenvalues $\lambda_{\rm e1}$ and $\lambda_{\rm e2}$ are always negative. Hence, we proved that the objective function of problem $\left(P4\right)$ is a concave function. Besides, it is noticed that constraints $C2$ and $C10$ are linear.

 Next, we set $f_i = \Delta P_{\rm L,min} - \eta_{\rm h} P_{\rm a} \left(1-\Gamma^2_{{\rm a}i}\right)$ corresponds to constraint $C11$. The second derivative of $f_i$ with respect to $\Gamma_{{\rm a}i}$ is $\frac{\partial^2 f_i}{\partial \Gamma^2_{{\rm a}i}} = 2 \eta_{\rm h} P_{\rm a} \geq 0$, which implies constraint $C11$ is convex. Since the objective function is concave, and the constraints $C2$, $C10$, and $C11$ are either linear or convex, we can conclude that problem $\left(P4\right)$ is a convex optimization problem. Hence, we proved Lemma \ref{L4}.

\vspace{-2mm}

\section{Maximization of Average Harvested Power in Nonlinear Energy Harvesting Models}\label{ApdC}

Here, we examine the harvested power of a passive tag using a saturation nonlinear model. The harvested power when backscattering symbol $S_i$ is given by~\cite{clerckx2019fundamentals}:
\begin{equation}
    P_{{\rm nL}i} = \frac{ \Psi_i - P_{\rm sat} \Omega }{1 - \Omega}, \quad i \in \mathbb{M}, 
    \label{eqA8}
\end{equation}
where $P_{\rm sat}$ is the maximum harvested power when the energy harvesting circuit is saturated due to exceedingly high input RF power, with 
\begin{align}
    \Omega &= \frac{1}{ 1 + e^{\varkappa \varpi} }, \label{eqA9} \\
    \Psi_i &= \frac{P_{\rm sat}}{ 1 + e^{ -\varkappa \left( P_{{\rm T}i} - \varpi  \right) } }, \quad i \in \mathbb{M}, \label{A10}
\end{align}
with $\varkappa$ and $\varpi$ being constants related to the circuit configuration and determined through standard curve fitting tools~\cite{boshkovska2015pracatical}. It is noteworthy that~\eqref{A10} is a sigmoid function influenced by the received RF power delivered to the rectifier, where $P_{{\rm T}i} = P_{\rm a} \left(1 - \Gamma^2_{{\rm nl},{\rm a}i} \right)$. We denote $\Gamma_{{\rm nl},{\rm a}i}$ as the reflection coefficient for the nonlinear model (similar to $\Gamma_{{\rm a}i}$ as in the linear model) to avoid confusion. 

Next, we formulate an optimization problem to maximize the average harvested power $P_{\rm nL,avg}$ by finding the optimal reflection coefficients. We consider the same operational requirements and constraints as in problem $\left(P2\right)$, except replacing $C8$ with $P_{{\rm nL}i} \geq P_{\rm L,min}$, expressed as follows:
\begin{align*}
(P7): &\max_{ \mathbf{\Gamma_{\rm nl}} } \quad P_{\rm nL,avg} = \smashoperator[r]{\sum_{i=1}^M} p_i P_{{\rm nL}i} \\
\text{s.t.} \hspace{3mm} &C16: \Gamma_{{\rm nl},{\rm a}i} \in \left[ 0,1 \right], \forall i \in \mathbb{M}, \\
    &C17: \frac{\abs{ \Gamma_{{\rm nl},{\rm a}i}-\Gamma_{{\rm nl},{\rm a}k} }}{2} \geq m_{\rm th},  \hspace{2mm} \forall i,k \in \mathbb{M}, \hspace{1mm} i \neq k,\\
    &C18: P_{{\rm nL}i} \geq P_{\rm L,min}, \hspace{2mm} \forall i \in \mathbb{M},\\
    &C19: \eta_{\rm b} P_{\rm a} G_{\rm r} \left( 1 - \Gamma_{{\rm nl},{\rm a}i} \right)^2 \geq P_{\rm b,min}, \hspace{2mm} \forall i \in \mathbb{M},
\end{align*}
where $\mathbf{\Gamma_{\rm nl}} = \left[ \Gamma_{{\rm nl},{\rm a1}}, \Gamma_{{\rm nl},{\rm a2}}, ..., \Gamma_{{\rm nl},{\rm a}M} \right]$. Remarkably, our proposed Algorithm~\ref{Algo1} to determine the optimal reflection coefficients that maximize the average harvested power applies to all energy harvest models that are increasing functions of the input power. Given $P_{\rm nL,avg}$ increases with $P_{{\rm T}i}$, we can use the solution methodology discussed in Section~\ref{Solution} to obtain the maximum harvested power $P^\ast_{\rm nL,avg}$. The key difference lies in the upper and lower bounds of the reflection coefficient. Specifically, we define $\Gamma_{\rm nl,ub}$ and $\Gamma_{\rm nl,lb}$ as the upper and lower bounds of $\Gamma_{{\rm nl},{\rm a}i}$, which can be determined from constraints $C16 - C19$. These bounds are given by:
\begin{align}
    \Gamma_{\rm nl,ub} &= \min \Biggl\{ \sqrt{\mathcal{X}} , 1 - \sqrt{\scalebox{1.1}{$\frac{P_{\rm b,min}}{ \eta_{\rm b} P_{\rm a} G_{\rm r}}$}} \Biggr\}, \\
    \Gamma_{\rm nl,lb} &= -\sqrt{\mathcal{X}},
\end{align}
where $\mathcal{X}$ is obtained by solving $\Gamma^2_{{\rm nl},{\rm a}i}$ from $P_{{\rm nl},{\rm L}i} = P_{\rm L,min}$ using numerical methods. Subsequently, we apply the technique introduced in Section~\ref{secSolution}. Specifically, we first determine $\Gamma_{{\rm nl},(i)}$ corresponding to the symbol with the occurrence probability $p_{(i)}$, similar to how $\Gamma_{{\rm nl},{\rm a}i}$ corresponds to $p_i$, with $\Gamma_{{\rm nl},\left(i\right)} > \Gamma_{{\rm nl},\left(i+1\right)}$ and $\Gamma_{{\rm nl},\left(i\right)} - \Gamma_{{\rm nl},\left(i+1\right)} = 2 m_{\rm th}$. Once we determine all the optimal $\Gamma_{{\rm nl},\left(i\right)}$ $\left(\right. \hspace{-1mm} \text{denoted as } \Gamma^\ast_{{\rm nl},\left(i\right)} \hspace{-1mm} \left.\right)$ for a particular probability sequence, we can then determine $\Gamma^\ast_{{\rm nl},{\rm a}i}$ by iteratively finding $\Gamma^\ast_{{\rm nl},\left(i\right)}$, $\forall i \in \mathbb{M}$ that results in the highest $P_{\rm nL,avg}$ among all probability sequences. Therefore, we express $P_{\rm nL,avg}$ in terms of $\Gamma_{{\rm nl},\left(1\right)}$, denoted as $P_{\rm \left(nL,avg\right)}$, as follows:
\begin{align}
    &P_{\rm \left(nL,avg\right)} \nonumber \\
    &= \frac{ \smashoperator[r]{\sum_{i=1}^M} \frac{p_{(i)} P_{\rm sat}}{ 1 + e^{ -\varkappa \left\{ P_{\rm a} \left[ 1 - \left( \Gamma_{{\rm nl},\left(1\right)} - 2 m_{\rm th} \left(i-1\right)  \right)^2 \right] - \varpi  \right\} } } - P_{\rm sat} \Omega }{1 - \Omega}. 
    \label{eqA13}
\end{align}

At this point, the original problem $\left(P7\right)$ can be transformed into a single-variable problem, and all constraints can be replaced by the boundary condition of the reflection coefficient. Given $P_{\rm sat}$ and $\Omega$ in~\eqref{eqA13} are constant, $\Gamma^\ast_{{\rm nl},\left(1\right)}$ can be obtained by solving the following problem:
\begin{align*}
(P8): &\operatorname*{arg\,max}_{\Gamma_{{\rm nl},\left(1\right)}} \smashoperator[r]{\sum_{i=1}^M} \frac{ p_{(i)} P_{\rm sat}}{ 1 + e^{ -\varkappa \left\{ P_{\rm a} \left[ 1 - \left( \Gamma_{{\rm nl},\left(1\right)} - 2 m_{\rm th} \left(i-1\right)  \right)^2 \right] - \varpi  \right\} } } \\
\text{s.t.} \hspace{3mm} &C20: \Gamma_{{\rm nl},\left(1\right)} - 2 m_{\rm th} \left(M-1\right) \geq \Gamma_{\rm lb}, \\
    &C21: \Gamma_{{\rm nl},\left(1\right)} \leq \Gamma_{\rm ub},
\end{align*}
where maximizing the objective function of problem $\left(P8\right)$ is equivalent to maximizing $P_{\rm \left(nL,avg\right)}$, and constraints $C20$ and $C21$ define the upper and lower bounds of $\Gamma_{{\rm nl},\left(1\right)}$. Given that the nonconvex problem $\left(P8\right)$ has only one variable, we employ a brute force method to determine $\Gamma^\ast_{{\rm nl},\left(1\right)}$ and then calculate $P^\ast_{\rm \left(nL,avg\right)}$, with the required time being insignificant. Next, we iteratively solve problem $\left(P8\right)$ for all possible probability sequences in each row of matrix $\mathbf{A}$. The probability sequence and the corresponding $\Gamma_{{\rm nl},\left(1\right)}$ that results in the highest $P^\ast_{\left(\rm nL,avg\right)}$ is the solution. In other words, Algorithm~\ref{Algo1} can be modified to solve problem $\left(P7\right)$ by replacing $\Gamma^\ast_{(i)}$, $P^\ast_{\rm \left(L,avg\right)}$, $\Gamma^\ast_{{\rm a}i}$, $P^\ast_{\rm L,avg}$ with $\Gamma_{{\rm nl,}(i)}$, $P^\ast_{\left(\rm nL,avg\right)}$, $\Gamma^\ast_{{\rm nl,a}i}$, $P^\ast_{\rm nL,avg}$ , and also replace step 6 with the solution of problem $\left(P8\right)$.

\vspace{-2mm}



\end{appendices}

\bibliographystyle{IEEEtran}
\bibliography{Biblio.bib}

\end{document}